\documentclass{elsart}
\usepackage{amsmath, amssymb}

\renewcommand{\phi}{\varphi}
\renewcommand{\epsilon}{\varepsilon}

\newcommand{\BB}{\mathbb}
\newcommand{\cal}{\mathcal}
\newcommand{\g}{\mathfrak}
\newcommand{\separate}{\vskip11pt}
\newcommand{\supp}{\operatorname{supp}}

\journal{the Journal of Functional Analysis}

\begin{document}

\begin{frontmatter}

\title{\bf A Localization Argument for Characters of Reductive Lie Groups}

\author{Matvei Libine}

\address{Department of Mathematics and Statistics,
University of Massachusetts,
710~North~Pleasant~Street, Amherst, MA 01003, USA}

\ead{matvei@math.umass.edu}

\tableofcontents

\begin{abstract}
This article provides a geometric bridge between two
entirely different character formulas for reductive Lie
groups and answers the question posed by W.~Schmid in
(1997, Deformation Theory and Symplectic Geometry,
{\bf 20}, 259-270).

A corresponding problem in the compact group setting was
solved by N.~Berline, E.~Getzler and M.~Vergne in
(1992, Heat Kernels and Dirac Operators)
by an application of the theory of equivariant forms
and particularly the fixed point integral localization formula.
This article (besides its representation-theoretical significance)
provides a whole family of examples where it is possible
to localize integrals to fixed points with respect to an
action of a {\em noncompact} group.
Moreover, a localization argument given here is not specific to
the particular setting considered in this article and can be
extended to a more general situation.

There is a broadly accessible article
(M.~Libine, 2002, math.RT/0208024)
which explains how the argument works in the $SL(2,\BB R)$ case,
where the key ideas are not obstructed by technical details and
where it becomes clear how it extends to the general case.
\end{abstract}

\begin{keyword}
integral character formula, fixed point character formula,
characteristic cycles of sheaves, equivariant forms,
fixed point integral localization formula
\end{keyword}

\end{frontmatter}

\separate

\begin{section}
{Acknowledgments}
\end{section}

I would like to express my deep gratitude to my thesis advisor
Wilfried Schmid for teaching me how to approach research problems,
suggesting this topic for my dissertation,
introducing me to representation theory,
his constructive criticism and his moral support.

\separate

\begin{section}
{Introduction}
\end{section}

For motivation, let us start with the case of a compact group.
Thus we consider a connected compact group $K$ and a
maximal torus $T \subset K$.
Let $\g k_{\BB R}$ and $\g t_{\BB R}$ denote the Lie algebras
of $K$ and $T$ respectively, and $\g k$, $\g t$ be their
complexified Lie algebras.
Let $\pi$ be a finite-dimensional representation of $K$.
We define a {\em character} on $K$ by
$\Theta_{\pi}(x) =_{\text{def}} \operatorname{tr} (\pi(x))$,
$x \in K$.
Any finite-dimensional representation is completely determined
(up to isomorphism) by its character.
We use the exponential map $\exp : \g k_{\BB R} \to K$ to define the
{\em character on the Lie algebra} of the representation $\pi$:
$$
\theta_{\pi}= (\det \exp_*)^{1/2} \exp^* \Theta_{\pi}.
$$
It is a smooth bounded function on $\g k_{\BB R}$.
Because $K$ is connected and compact, the exponential map is
surjective and generically non-singular.
Thus $\theta_{\pi}$ still determines the representation.

Now let us assume that the representation $\pi$ is irreducible.
Then there are two entirely different character formulas
for $\theta_{\pi}$ -- the {\em Weyl character formula}
and {\em Kirillov's character formula}.
Recall that the irreducible representations of $K$
can be enumerated by their {\em highest weights} which are
elements of the {\em weight lattice} $\Lambda$ in $i \g t_{\BB R}^*$
intersected with a chosen {\em Weyl chamber}.
Let $\lambda=\lambda(\pi) \in i \g t_{\BB R}^*$
denote the highest weight corresponding to $\pi$.
Let $W=N_K(\g t_{\BB R})/T$,
where $N_K(\g t_{\BB R})$ is the normalizer of $\g t_{\BB R}$ in $K$.
The set $W$ is a finite group called the {\em Weyl group};
it acts on $\g t_{\BB R}$ and hence on $i \g t_{\BB R}^*$.
We can choose a positive definite inner product
$\langle \cdot,\cdot\rangle$ on $i \g t_{\BB R}^*$ invariant under $W$.
Then the Weyl character formula can be stated as follows:
$$
\theta_{\pi}|_{\g t_{\BB R}}(t)=\sum_{w \in W} \frac {e^{w(\lambda+\rho)(t)}}
{\prod_{\alpha \in \Phi, \langle w(\lambda+\rho),\alpha\rangle>0} \alpha(t)},
$$
where $\Phi \subset i \g t_{\BB R}^*$ is the {\em root system}
of $\g k_{\BB R}$, and
$\rho \in i \g t_{\BB R}^*$ is a certain small vector independent
of $\pi$. Because $\theta_{\pi}$ is $Ad(K)$-invariant and
every $Ad(K)$-orbit in $\g k_{\BB R}$ meets $\g t_{\BB R}$,
this formula completely determines $\theta_{\pi}$.

Kirillov's character formula provides a totally different expression
for the irreducible characters on $\g k_{\BB R}$.
The splitting
$\g k_{\BB R}=\g t_{\BB R} \oplus [\g t_{\BB R},\g k_{\BB R}]$
induces a dual splitting of the vector space $i \g k_{\BB R}^*$,
which allows us to think of $\lambda$ and $\rho$ as lying in
$i \g k_{\BB R}^*$. The adjoint action of $K$ on
$\g k_{\BB R}$ has a dual action on $i \g k_{\BB R}^*$
called coadjoint representation. We define
$$
\Omega_{\lambda+\rho} = \text{ $K$-orbit of $\lambda+\rho$
in $i \g k_{\BB R}^*$}.
$$
It will be convenient to define the Fourier transform $\hat \phi$
of a test function $\phi \in {\cal C}^{\infty}_c (\g k_{\BB R})$
without the customary factor of $i=\sqrt{-1}$ in the exponent,
as a function on $i \g k_{\BB R}^*$:
$$
\hat \phi(\zeta) = \int_{\g k_{\BB R}}
\phi(x)e^{\langle \zeta, x \rangle} dx.
$$
Then Kirillov's character formula describes $\theta_{\pi}$ as a
distribution on $\g k_{\BB R}$:
$$
\int_{\g k_{\BB R}} \theta_{\pi} \phi dx =
\int_{\Omega_{\lambda+\rho}} \hat \phi d\beta,
$$
where $d\beta$ is the measure induced by the canonical symplectic
structure of $\Omega_{\lambda+\rho}$.
In other words,
$$
\hat \theta_{\pi} = \text{ integration over $\Omega_{\lambda+\rho}$}.
$$
Kirillov calls this the ``universal formula'' for irreducible
characters.

The geometric relationship between these two formulas is even
more striking. As a homogeneous space, $\Omega_{\lambda+\rho}$
is isomorphic to the {\em flag variety} $X$, i.e the variety
of Borel subalgebras
$\g b \subset \g k = \g k_{\BB R} \otimes_{\BB R} \BB C$.
The Borel-Weil-Bott theorem can be regarded as an explicit
construction of a holomorphic $K$-equivariant line
bundle ${\cal L}_{\lambda} \to X$ such that the resulting
representation of $K$ in the cohomology groups is:
\begin{eqnarray*}
H^p(X, {\cal O}({\cal L}_{\lambda})) &=& 0 \text{\quad if $p \ne 0$},  \\
H^0(X, {\cal O}({\cal L}_{\lambda})) &\simeq& \pi.
\end{eqnarray*}
Then the Weyl character formula is a consequence of
the Atiyah-Bott fixed point formula.
On the other hand, N.~Berline, E.~Getzler
and M.~Vergne proved in \cite{BGV}
Kirillov's character formula using the integral
localization formula for $K$-equivariant forms.
They showed that the right hand side of Kirillov's character formula
equals the right hand side of the Weyl character formula.
Recall that integral of an equivariant
form is a function on $\g k_{\BB R}$, and
the localization formula
reduces integration of an equivariantly closed form to
summation over the zeroes of the vector field in $X$ generated
by $k \in \g k_{\BB R}$. It is crucial for the localization
formula to hold that the group $K$ is compact.

\separate

Now let $G_{\BB R}$ be a connected, linear, reductive Lie group.
We let $\g g_{\BB R}$ denote its Lie algebra.
Then most representations of interest have infinite dimension.
We always consider representations on complete, locally convex
Hausdorff topological vector spaces and require that the action
of $G_{\BB R}$ is continuous.
Let $K$ be a maximal compact subgroup of $G_{\BB R}$.
A reasonable category of representations consists of
{\em admissible} representations of {\em finite length}.
(A representation $\pi$ has finite length if every increasing chain
of closed, invariant subspaces breaks off after finitely many steps;
$\pi$ is admissible if its restriction to $K$ contains any irreducible
representation of $K$ at most finitely often.) Admissibility
is automatic for irreducible unitary representations.
Although trace of a linear operator in an infinite-dimensional
space cannot be defined in general, it is still possible to
define a character $\theta_{\pi}$ as an $Ad(G_{\BB R})$-invariant
distribution on $\g g_{\BB R}$. (See \cite{A} for details.)

M.~Kashiwara and W.~Schmid in their paper \cite{KSch}
generalize the Borel-Weil-Bott construction.
Instead of line bundles on the flag variety $X$ they consider
$G_{\BB R}$-equivariant sheaves ${\cal F}$ and, for each integer
$p \in \BB Z$, they define representations of $G_{\BB R}$ in
$\operatorname{Ext}^p({\cal F},{\cal O})$.
Such representations turn out to be admissible of finite length.
Then W.~Schmid and K.~Vilonen prove in \cite{SchV2}
two character formulas for these representations --
the fixed point character formula and the integral character formula.
In the case when $G_{\BB R}$ is compact, the former reduces to
the Weyl character formula and the latter -- to Kirillov's character
formula. The fixed point formula was conjectured by M.~Kashiwara
in \cite{K}, and its proof uses a
generalization of the Lefschetz fixed point formula to sheaf
cohomology due to M.~Goresky and R.~MacPherson in \cite{GM}.
On the other hand, W.~Rossmann in \cite{R} established existence
of an integral character formula over an unspecified
Borel-Moore cycle.
W.~Schmid and K.~Vilonen prove the integral character formula
where integration takes place over the characteristic cycle of
${\cal F}$, $Ch({\cal F})$, and their proof
depends totally on representation theory.

\separate

The equivalence of these two formulas can be stated in terms
of the sheaves ${\cal F}$ alone, without any reference to
their representation-theoretic significance.
In the announcement \cite{Sch} W.~Schmid posed a question:
``Can this equivalence be seen directly without a detour to
representation theory, just as in the compact case.''

In this article I provide such a geometric link.
I introduce a localization technique which allows
to localize integrals to the zeroes of vector fields on $X$
generated by the infinitesimal action of $\g g_{\BB R}$.
Thus, in addition to a representation-theoretical result,
we obtain a whole family of examples where it is possible
to localize integrals to fixed points with respect to an
action of a noncompact group.
Moreover, the localization argument given here is not specific to
this particular setting and can be extended to a more
general situation.

There is a broadly accessible article \cite{L}
which explains how the argument works in the $SL(2,\BB R)$ case,
where the key ideas are not obstructed by technical details and
where it becomes clear how it extends to the general case.

\separate

\begin{section}
{Setup}
\end{section}

In these notes we try to keep the same notations as W.~Schmid and K.~Vilonen
use in \cite{SchV2} as much as possible.
That is they fix a connected, complex algebraic,
reductive group $G$ which is defined over $\BB R$.
The representations they consider are representations of a real form
$G_{\BB R}$ of $G$ -- in other words, $G_{\BB R}$ is a subgroup of
$G$ lying between the group of real points $G(\BB R)$ and the
identity component $G(\BB R)^0$.
They regard $G_{\BB R}$ as a reductive Lie group and
denote by $\g g$ and $\g g_{\BB R}$ the Lie algebras of
$G$ and $G_{\BB R}$ respectively,
they also denote by $X$ the flag variety of $G$.

If $g \in \g g$ is an element of the Lie algebra,
we denote by $\operatorname{VF}_g$ the vector field on $X$
generated by $g$: if $x \in X$ and $f \in {\cal C}^{\infty}(X)$, then
$$
\operatorname{VF}_g (x)f=
\frac d{d\epsilon} f(\exp(\epsilon g) \cdot x)|_{\epsilon=0}.
$$
We call a point $x \in X$ a {\em fixed point of $g$} if the vector field
$\operatorname{VF}_g$ on $X$ vanishes at $x$, i.e.
$\operatorname{VF}_g (x)= 0$.

In this paragraph we explain the general picture, but since objects
mentioned here will not play any role in what follows they will not
be defined, rather the reader is referred to \cite{SchV2}.
W.~Schmid and K.~Vilonen denote by $\g h$ the universal Cartan
algebra. They pick an element $\lambda \in \g h^*$ and introduce
the ``$G_{\BB R}$-equivariant derived category on $X$ with twist
$(\lambda - \rho)$'' denoted by
$\operatorname{D}_{G_{\BB R}}(X)_{\lambda}$.
They also introduce ${\cal O}_X(\lambda)$,
the twisted sheaf of holomorphic functions on $X$,
with twist $(\lambda - \rho)$.
Then, for ${\cal F} \in \operatorname{D}_{G_{\BB R}}(X)_{-\lambda}$,
they define a virtual representation of $G_{\BB R}$
\begin{equation}  \label{rep}
\sum_p (-1)^p \operatorname{Ext}^p (\BB D {\cal F}, {\cal O}_X(\lambda)),
\end{equation}
where $\BB D {\cal F} \in \operatorname{D}_{G_{\BB R}}(X)_{\lambda}$
denotes the Verdier dual of ${\cal F}$.
It was shown in \cite{KSch} that each
$\operatorname{Ext}^p (\BB D {\cal F}, {\cal O}_X(\lambda))$
is admissible of finite length.
In particular, this representation has a $\g g_{\BB R}$-character $\theta$.
We think of a character as an $Ad(G_{\BB R})$-invariant linear functional
defined on the space of smooth compactly supported differential forms
$\phi$ on $\g g_{\BB R}$ of top degree, and write
$\int_{\g g_{\BB R}} \theta \phi$ for the value of $\theta$ at $\phi$.

Then Proposition 3.1 and Theorem 3.3 of \cite{A} say that
the character $\theta$ is given by integration against a function
$F_{\theta} \in L^1_{loc}(\g g_{\BB R})$:
$$
\int_{\g g_{\BB R}} \theta \phi = \int_{\g g_{\BB R}} F_{\theta} \phi.
$$
This function $F_{\theta}$ is invariant under the adjoint action of
$G_{\BB R}$ on $\g g_{\BB R}$, and the restriction of $F_{\theta}$ to
the set of regular semisimple elements of $\g g_{\BB R}$ can be
represented by an analytic function.

There are two formulas expressing the character $\theta$ of the
virtual representation (\ref{rep}) as a distribution on $\g g_{\BB R}$.
We will start with the right hand side of the integral character formula
$$
\int_{\g g_{\BB R}} \theta \phi =
\frac 1{(2\pi i)^nn!} \int_{Ch({\cal F})}
\mu_{\lambda}^* \hat \phi (-\sigma+\pi^* \tau_{\lambda})^n
$$
and show that it is equivalent to the right hand side of the
fixed point character formula
$$
\int_{\g g_{\BB R}} \theta \phi = \int_{\g g_{\BB R}} F_{\theta} \phi,
\qquad
F_{\theta} (g) =
\sum_{k=1}^{|W|} 
\frac {m_{x_k(g)} e^{\langle g,\lambda_{x_k(g)} \rangle}}
{\alpha_{x_k(g),1}(g) \dots \alpha_{x_k(g),n}(g)},
$$
where $x_1, \dots, x_k$ are the fixed points of $g$ and
integers $m_{x_k(g)}$'s are the local invariants of ${\cal F}$.
This way we will obtain a new proof of the integral character formula.

\separate

Because both character formulas depend on
${\cal F} \in \operatorname{D}_{G_{\BB R}}(X)_{-\lambda}$ only through
its characteristic cycle $Ch({\cal F})$, we can simply
replace ${\cal F}$ with a $G_{\BB R}$-equivariant sheaf on the flag
variety $X$ with the same characteristic cycle.
We will use the same notation ${\cal F}$ to denote this
$G_{\BB R}$-equivariant sheaf on $X$.
Let $n = \dim_{\BB C} X$, let $\pi: T^*X \twoheadrightarrow X$ be
the projection map, and equip $\g g_{\BB R}$ with some orientation.
We will make an elementary calculation of the integral
\begin{equation}   \label{integral_formula}
\frac 1{(2\pi i)^nn!} \int_{Ch({\cal F})}
\mu_{\lambda}^* \hat \phi (-\sigma+\pi^* \tau_{\lambda})^n
\end{equation}
where $\phi$ is a smooth compactly supported differential
form on $\g g_{\BB R}$ of top degree,
$$
\hat \phi(\zeta) = \int_{\g g_{\BB R}} e^{\langle g, \zeta \rangle} \phi
\qquad (g \in \g g_{\BB R},\,\zeta \in \g g^*)
$$
is its Fourier transform
(without the customary factor of $i=\sqrt{-1}$ in the exponent),
$\mu_{\lambda}: T^*X \to \g g^*$
is the twisted moment map defined in \cite{SchV1}
and $\tau_{\lambda}$, $\sigma$ are 2-forms on
$X$ and $T^*X$ respectively defined in \cite{SchV2}.
The form $\sigma$ is the complex algebraic symplectic form on $T^*X$.
On the other hand, the precise definition of the form $\tau_{\lambda}$
will not be important. What will be important, however,
is that, for each $g \in \g g$, the $2n$-form on $T^*X$
\begin{equation}  \label{integrand}
e^{\langle g, \mu_{\lambda}(\zeta) \rangle}
(-\sigma+\pi^* \tau_{\lambda})^n
\end{equation}
is closed.
By Lemma 3.16 in \cite{SchV2}, the integral (\ref{integral_formula})
converges absolutely. This is true essentially because the Fourier
transform $\hat \phi$ is holomorphic on $\g g^*$ and decays rapidly in
the imaginary directions.
Because the form (\ref{integrand}) is not $Ad(G_{\BB R})$-invariant,
at this point it is not even clear that the distribution defined by
(\ref{integral_formula}) is $Ad(G_{\BB R})$-invariant.

\begin{rem}
One can give a direct proof that the distribution
$$
\tilde \theta: \phi \mapsto \tilde \theta(\phi) =
\frac 1{(2\pi i)^nn!} \int_{Ch({\cal F})}
\mu_{\lambda}^* \hat \phi (-\sigma+\pi^* \tau_{\lambda})^n
$$
is $Ad(G_{\BB R})$-invariant and that
$\tilde \theta$ is an eigendistribution of every biinvariant differential
operator on $G_{\BB R}$. Then  Theorem 3.3 from \cite{A} will imply
that there exists a function $\tilde F \in L_{loc}^1 (\g g_{\BB R})$
such that
$$
\tilde \theta (\phi) = \int_{\g g_{\BB R}} \tilde F \phi,
$$
this function $\tilde F$ is invariant under the adjoint action of
$G_{\BB R}$ on $\g g_{\BB R}$, and the restriction of $\tilde F$ to
the set of regular semisimple elements of $\g g_{\BB R}$ can be
represented by an analytic function.

However, this information would not make our computation of integral
(\ref{integral_formula}) any easier and it will follow automatically
from the corresponding statement for the fixed point character formula.
\end{rem}

\separate

When $\lambda \in \g h^*$ is regular, the twisted moment map
$\mu_{\lambda}$ is a real analytic
diffeomorphism of $T^*X$ onto $\Omega_{\lambda} \subset \g g^*$
-- the orbit of $\lambda$ under the coadjoint action of $G$ on $\g g^*$.
Let $\sigma_{\lambda}$ denote the canonical $G$-invariant complex algebraic
symplectic form on $\Omega_{\lambda}$. Then
$e^{\langle g, \zeta \rangle} (\sigma_{\lambda})^n$ is a holomorphic
$2n$-form of maximal possible degree, hence closed.
Proposition 3.3 in \cite{SchV2} says that
$\mu_{\lambda}^* (\sigma_{\lambda}) = -\sigma+\pi^* \tau_{\lambda}$.
This shows that, for $\lambda$ regular, $g \in \g g$,
$$
e^{\langle g, \mu_{\lambda}(\zeta) \rangle}
(-\sigma+\pi^* \tau_{\lambda})^n =
\mu_{\lambda}^* \bigl( e^{\langle g, \zeta \rangle} (\sigma_{\lambda})^n \bigr)
$$
is a closed $2n$-form on $T^*X$.
Note that neither map $\mu_{\lambda}$ nor the form (\ref{integrand})
is holomorphic. Because the form
$e^{\langle g, \mu_{\lambda}(\zeta) \rangle}
(-\sigma+\pi^* \tau_{\lambda})^n$
depends on $\lambda$ real analytically and the set of regular elements
is dense in $\g h^*$, we conclude that the form in the equation
(\ref{integrand}) is closed.

If $\lambda \in \g h^*$ is regular, then we can rewrite our
integral (\ref{integral_formula}) as
$$
\frac 1{(2\pi i)^nn!} \int_{Ch({\cal F})}
\mu_{\lambda}^* \hat \phi (-\sigma+\pi^* \tau_{\lambda})^n =
\frac 1{(2\pi i)^nn!} \int_{(\mu_{\lambda})_* Ch({\cal F})}
\hat \phi (\sigma_{\lambda})^n.
$$
This is a Rossmann type character formula.
W.~Rossmann in \cite{R} established existence
of an integral character formula over an unspecified
Borel-Moore cycle in the coadjoint orbit $\Omega_{\lambda}$.
This expression tells us that we can choose Rossmann's cycle to be
$(\mu_{\lambda})_* Ch({\cal F})$.

\separate

Fix a norm $\|.\|$ on $\g g^*$. Then the moment map $\mu$
induces a vector bundle norm on $T^*X$: for $\zeta \in T^*X$
its norm will be $\|\mu(\zeta)\|$. We will use the same notation
$\|.\|$ for this norm too.

Let $\g g_{\BB R}'$ denote the set of regular semisimple elements
$g \in \g g_{\BB R}$ which satisfy the following additional property:
If $\g t_{\BB R} \subset \g g_{\BB R}$ and $\g t \subset \g g$ are
the unique Cartan subalgebras in $\g g_{\BB R}$ and $\g g$ respectively
containing $g$, $\alpha \in \g t^*$ is a (complex) root such that
$Re(\alpha)|_{\g t_{\BB R}} \not \equiv 0$, then $Re(\alpha(g)) \ne 0$.

Since the complement of $\g g_{\BB R}'$ in $\g g_{\BB R}$ has measure zero,
we can replace integration over $\g g_{\BB R}$ by
integration over $\g g_{\BB R}'$. Then
\begin{multline}  \label{theintegral}
\frac 1{(2\pi i)^nn!} \int_{Ch({\cal F})}
\mu_{\lambda}^* \hat \phi (-\sigma+\pi^* \tau_{\lambda})^n \\
=\lim_{R \to \infty}
\frac 1{(2\pi i)^nn!}
\int_{\g g_{\BB R} \times (Ch({\cal F}) \cap \{\|\zeta\| \le R\})}
e^{\langle g, \mu_{\lambda}(\zeta) \rangle} \phi
(-\sigma+\pi^* \tau_{\lambda})^n \\
=\lim_{R \to \infty} \frac 1{(2\pi i)^nn!}
\int_{\g g_{\BB R}' \times (Ch({\cal F}) \cap \{\|\zeta\| \le R\})}
e^{\langle g, \mu_{\lambda}(\zeta) \rangle} \phi
(-\sigma+\pi^* \tau_{\lambda})^n.
\end{multline}
(Of course, the orientation on
$\g g_{\BB R}' \times (Ch({\cal F}) \cap \{\|\zeta\| \le R\})$
is induced by the product orientation on
$\g g_{\BB R} \times Ch({\cal F})$.)

We will interchange the order of integration:
integrate over the characteristic cycle first and only then
perform integration over $\g g_{\BB R}'$.
The integrand is an equivariant form with respect to some
compact real form $U_{\BB R} \subset G$.
But $U_{\BB R}$ does not preserve $Ch({\cal F})$
(unless $Ch({\cal F})$ is a multiple of the zero section of $T^*X$
equipped with some orientation).
Each $g \in \g g_{\BB R}'$ has exactly $|W|$ fixed points on $X$,
where $|W|$ is the order of the Weyl group.
We regard the integral (\ref{theintegral})
as an integral of a closed differential form
$e^{\langle g, \mu_{\lambda}(\zeta) \rangle}
\phi(-\sigma+\pi^* \tau_{\lambda})^n$
over a chain in $\g g_{\BB R}' \times (T^*X \cap \{\|\zeta\| \le R\})$.

Because the closure of a $G_{\BB R}$-orbit on $X$ may be extremely
singular, same is true of $Ch({\cal F})$.
We will use the open embedding theorem of W.~Schmid and K.~Vilonen
(\cite{SchV1}) to construct a deformation of
$Ch({\cal F})$ into a simple cycle of the following kind:
$$
m_1 T^*_{x_1}X +\dots+ m_{|W|} T^*_{x_{|W|}}X,
$$
where $m_1,\dots,m_{|W|}$ are some integers,
$x_1,\dots,x_{|W|}$ are the points in $X$ fixed by $g \in \g g_{\BB R}'$,
and each cotangent space $T^*_{x_k}X$ is given some orientation.
This is very similar to the classical Morse's lemma which says
that if we have a smooth real valued function $f$ on a manifold
$M$, then the sublevel sets $\{ m \in M;\,f(m) < a \}$ and
$\{ m \in M;\,f(m) < b \}$ can be deformed one into the other
as long as there are no critical values of $f$ in an open
interval containing $a$ and $b$.
To ensure that the integral (\ref{theintegral})
behaves well, we will stay during the process of deformation inside the set
$$
\{(g,\zeta) \in \g g_{\BB R}' \times T^*X;\,
Re( \langle g,\mu(\zeta) \rangle) \le 0\}.
$$

After that we define another deformation $\Theta_t(g): T^*X \to T^*X$,
where $g \in \g g_{\BB R}'$, $t \in [0,1]$.
It has the following meaning. In the classical proof of
the Fourier inversion formula
$$
\phi(g)= \frac 1{(2\pi i)^{\dim_{\BB C}\g g}}
\int_{\zeta \in i \g g_{\BB R}^*}
\hat \phi(\zeta) e^{-\langle g,\zeta \rangle}
$$
we multiply the integrand by a term like $e^{-t \|\zeta\|^2}$
to make it integrable over $\g g_{\BB R} \times i \g g_{\BB R}^*$,
and then let $t \to 0^+$.
The deformation $\Theta_t(g)$ has a very similar effect --
it makes our integrand an $L^1$-object.
Lemma \ref{slanting} says that this substitute is permissible.
Its proof is very technical, but the idea is quite simple.
The difference between the original integral (\ref{theintegral})
and the deformed one is expressed by an integral of 
$e^{\langle g, \mu_{\lambda}(\zeta) \rangle}
\phi(-\sigma+\pi^* \tau_{\lambda})^n$ over a certain cycle
$\tilde C(R)$ supported in
$\g g_{\BB R}' \times (T^*X \cap \{\|\zeta\| = R\})$
which depends on $R$ by scaling along the fiber.
Recall that the Fourier transform $\hat \phi$ decays rapidly
in the imaginary directions which is shown by an integration
by parts. We modify this integration by parts argument to prove
a similar statement about behavior of the integrand on
the support of $\tilde C(R)$ as $R \to \infty$.
Hence the difference of integrals in question tends to zero.

The key ideas are the deformation of $Ch({\cal F})$,
the definition of $\Theta_t(g): T^*X \to T^*X$ and
Lemma \ref{slanting}. Because of the right definition
of $\Theta_t(g)$, Lemma \ref{slanting} holds and our calculation
of the integral (\ref{theintegral}) becomes very simple.
We will see that, as $R \to \infty$ and $t \to 0^+$,
the integral will concentrate more and more inside $T^*U$,
where $U$ is a neighborhood of the set of fixed points of $g$ in $X$.
In the limit, we obtain the right hand side of the
fixed point character formula.
This means that the integral (\ref{integral_formula})
is localized at the fixed points of $g$.

\separate

The following convention will be in force throughout these notes:
whenever $A$ is a subset of $B$,
we will denote the inclusion map $A \hookrightarrow B$ by
$j_{A \hookrightarrow B}$.

\separate

\begin{section}
{Deformations of Characteristic Cycles}
\end{section}

Recall the notion of families of cycles introduced by W.~Schmid and
K.~Vilonen in \cite{SchV1}.

Let $M$ be a complex manifold, ${\cal G}$ a constructible sheaf on
$M$ and $U$ a constructible open subset of $M$.
Then we denote by ${\cal G}_U$ the sheaf
$$
(j_{U \hookrightarrow M})_! \circ (j_{U \hookrightarrow M})^* ({\cal G}).
$$
We can calculate the characteristic cycle of ${\cal G}_U$
using the open embedding theorem of W.~Schmid and K.~Vilonen (\cite{SchV1}).
So let $f$ be a constructible function defining $U$.
That is $f$ is a constructible real-valued ${\cal C}^2$-function
defined on an open neighborhood of the closure $\overline{U} \subset M$
such that $f$ is strictly positive on $U$ and the boundary
$\partial U$ is precisely the zero set of $f$.
We regard $df$ as a subset of $T^*M$.
Suppose $a>0$ is such that $|Ch({\cal G})|$ and $df$
do not intersect over $V= \{ u \in U; 0< f(u) < a \}$, i.e.
$|Ch({\cal G})| \cap df \cap T^*V = \varnothing$.
For $\epsilon \in [0,a)$ define
$U_{\epsilon}=\{u \in U; \, f(u)>\epsilon\}$; in particular $U_0=U$.
Now, for each $\epsilon \in (0,a)$, we can consider a cycle
$Ch({\cal G}_{U_{\epsilon}})$ in $T^*M$.
Next proposition asserts that these cycles piece together to
form a family of cycles in $T^*M$ and this family has limit
$Ch({\cal G}_U)$ as $\epsilon \to 0^+$.

\separate

\begin{prop} \label{deformation1}
Suppose that $|Ch({\cal G})| \cap df \cap T^*V = \varnothing$.
Then there exists a family $C_{(0,a)}$ of ($\dim_{\BB R}M$)-dimensional
cycles in $T^*M$
parameterized by $(0,a)$ such that, for each $\epsilon \in (0,a)$,
the specialization at $\epsilon$, $C_{\epsilon}$, is equal to
$Ch({\cal G}_{U_{\epsilon}})$.
Moreover,
$$
\lim_{\epsilon \to 0^+} C_{\epsilon}=
\lim_{\epsilon \to 0^+} Ch({\cal G}_{U_{\epsilon}})=
Ch({\cal G}_U).
$$
\end{prop}

\pf
Let $p: (0,a) \times M \twoheadrightarrow M$ be the projection.
We have a function $\tilde f(\epsilon,m)=f(m)-\epsilon$
defined on an open neighborhood of the closure
$\overline{(0,a) \times U}$ in $(0,a) \times M$.
It defines an open set
$$
\tilde U=\{ x \in (0,a) \times U;\, \tilde f(x)>0\}
=\{(\epsilon,u) \in (0,a) \times U;\, f(u)>\epsilon\}.
$$
Let us consider a family of cycles $C''_I$ in $T^*(0,a) \times T^*M$
parameterized by some open interval $I=(0,b)$
calculating the cycle $Ch((p^*{\cal G})_{\tilde U})$:
$$
C''_s=Ch((p^*{\cal G})|_{\tilde U})-s \frac{d\tilde f}{\tilde f}, \qquad
\lim_{s \to 0^+} C''_s= Ch((p^*{\cal G})_{\tilde U}).
$$

\begin{lem} \label{proper}
The projection map
$\tilde p:
I \times T^*(0,a) \times T^*M \twoheadrightarrow I \times (0,a) \times T^*M$
is proper on the support $|C''_I|$.
\end{lem}

This is where the assumption that the intersection
$|Ch({\cal G})| \cap df \cap T^*V$ is empty is used.
We will assume this lemma for now and give an argument later.

Because the projection map
$\tilde p:
I \times T^*(0,a) \times T^*M \twoheadrightarrow I \times (0,a) \times T^*M$
is proper on the support $|C''_I|$, $\tilde p$ sends $C''_I$ into a cycle
$C'_I$ in $I \times (0,a) \times T^*M$ which is a family of cycles
in $(0,a) \times T^*M$.
Let $C_{(0,a)}=\lim_{s \to 0^+} C'_I$.
It is clear that $C_{\epsilon}=Ch({\cal G}_{U_{\epsilon}})$
which proves the first part of the proposition.

To prove that $\lim_{\epsilon \to 0^+} C_{\epsilon}=Ch({\cal G}_U)$
we need to show that if we regard $C_{(0,a)}$ as a chain in
$[0,a) \times T^*M$, then $\partial C_{(0,a)}=-Ch({\cal G}_U)$.
We can also view $C'_I$ as a family of cycles
in $I \times T^*M$ parameterized by $(0,a)$.
Let $C_I=\lim_{\epsilon \to 0^+} C'_I$.
It is clear that $C_I$ is nothing else but
$- \bigl( Ch({\cal G}|_U) - s \frac {df}f \bigr) $.
(The negative sign appears because of orientation matters.)
If we consider $C'_I$ as a chain in $[0,b) \times [0,a) \times T^*M$, then
its boundary $\partial C'_I$ is a cycle which can be written as a sum
of two chains: $\partial C'_I=-C_{(0,a)}-C_I$.
This is a cycle and there is no boundary, hence:
$$
\lim_{\epsilon \to 0^+} C_{\epsilon}=
-\partial C_{(0,a)}= \partial C_I =-\lim_{s \to 0^+} C_I
= \lim_{s \to 0^+} \bigl( Ch({\cal G}|_U) - s \frac {df}f \bigr) =
Ch({\cal G}_U).
$$

It remains to prove Lemma \ref{proper}.

\noindent{\it Proof of Lemma \ref{proper}.}
Fix a smooth inner product on $T^*M$.
It induces a metric which we denote by $\|.\|$.
To prove that the projection map
$\tilde p:
I \times T^*(0,a) \times T^*M \twoheadrightarrow I \times (0,a) \times T^*M$
is proper on the support $|C''_I|$
we need to show that the inverse image of a compact set
$\tilde K \subset I \times (0,a) \times T^*M$ is compact.
We can always enlarge $\tilde K$ which allows us to assume that
$$
\tilde K=
[s_0, s_1] \times [\epsilon_0, \epsilon_1] \times
\{\zeta;\,\zeta \in T^*_mM, m \in K, \|\zeta\| \le R \}
$$
for some $0<s_0<s_1<b$, $0<\epsilon_0<\epsilon_1<a$, $R>0$
and some compact set $K \subset M$.

There exists an angle $\alpha>0$ such that
whenever $m \in K \cap \operatorname{supp}{\cal G} \cap U$ and
$\epsilon_0 \le f(m) \le \frac{\epsilon_1+a}2$,
the open cone around $df(m)$ in $T^*_mM$,
\begin{multline*}
\operatorname{Cone}_{\alpha}(df(m))  \\
= \{ \zeta \in T^*_mM \setminus \{0\};
\text{ the angle between $\zeta$ and $df(m)$ is less than $\alpha$} \},
\end{multline*}
does not intersect $|Ch({\cal G})|$.
We will assume that $\alpha < \pi/2$.

The preimage of $\tilde K$ is
\begin{multline*}
\bigl\{ (r d\epsilon,s,\epsilon,\zeta) \in \BB R d\epsilon \times \tilde K;\,
\zeta \in T^*U,
(r d\epsilon,\epsilon,\zeta) \in \{ |Ch((p^*{\cal G})|_{\tilde U})|
- s\frac{d \tilde f}{\tilde f} \} \bigr\} \\
= \bigl\{ (r d\epsilon,s,\epsilon,\zeta) \in \BB R d\epsilon \times \tilde K;
\hspace{3in} \\
\zeta \in T^*U,
(r d\epsilon,\epsilon,\zeta) \in \{ |Ch((p^*{\cal G})|_{\tilde U})|
- s\frac{df - d\epsilon}{f-\epsilon} \} \bigr\} \\
= \bigl\{ (r d\epsilon,s,\epsilon,\zeta) \in \BB R d\epsilon \times \tilde K;
\hspace{3in} \\
\zeta \in T^*_mM, f(m)> \epsilon,
\zeta \in \{ |Ch({\cal G}|_U)|- s\frac{df}{f-\epsilon} \},
r = \frac s{f-\epsilon} \bigr\}.
\end{multline*}
When $f(m) > \frac{\epsilon_1+a}2$,
$r$ is bounded from above by $s_1 \frac 2{a-\epsilon_1}$.
Suppose now that $f(m) \le \frac{\epsilon_1+a}2$.
Conditions $\zeta \in T^*_mM$,
$\zeta \in \{ |Ch({\cal G}|_U)|- s\frac{df}{f-\epsilon} \}$ and
$|Ch({\cal G})| \cap \operatorname{Cone}_{\alpha}(df(m)) = \varnothing$
mean that $\zeta = \tilde \zeta - s\frac{df}{f-\epsilon}$ for some
$\tilde \zeta \in |Ch({\cal G}|_U)|$ and the angle between
$\tilde \zeta$ and $s\frac{df}{f-\epsilon}$ is at least $\alpha$.
It follows that
$$
\sin \alpha \Bigl\| s\frac{df}{f-\epsilon}(m) \Bigr\| \le \|\zeta\| \le R,
$$
which in turn implies
$$
s \frac{\|df\|}{f-\epsilon}(m) \le \frac R{\sin \alpha}.
$$

Let $D$ be the minimum of $\|df\|$ on
$$
K \cap \operatorname{supp} ({\cal G}) \cap
\bigl\{ m \in U;\, \epsilon_0 \le f(m) \le \frac{\epsilon_1+a}2 \bigr\},
$$
$D>0$.
Then $r=\frac s{f-\epsilon}$ can be at most 
$\frac R{D \sin \alpha}$.
On the other hand, $r \ge 0$.
Thus $r$ lies between $0$ and
$\max \bigl\{ s_1 \frac 2{a-\epsilon_1},\frac R{D \sin \alpha} \bigr\}$.
This proves that the preimage of
$\tilde K$ has bounded $r$-coordinate, hence compact.
This finishes the proof of Lemma \ref{proper} and the proposition.
\qed

For each $\epsilon \in (0,a)$, we can consider a different cycle:
$$
Ch \bigl( (Rj_{U_\epsilon \hookrightarrow M})_* \circ
(j_{U_\epsilon \hookrightarrow M})^* ({\cal G}) \bigr) =
Ch(R\Gamma_{U_{\epsilon}}{\cal G}).
$$
There is a similar result about these cycles too which we will also use.
Notice that in the intersection condition the section
$df$ is replaced with $-df$.

\begin{prop} \label{deformation2}
Suppose that $|Ch({\cal G})| \cap -df \cap T^*V = \varnothing$.
Then there exists a family $C_{(0,a)}$ of ($\dim_{\BB R}M$)-dimensional
cycles in $T^*M$
parameterized by $(0,a)$ such that, for each $\epsilon \in (0,a)$,
the specialization at $\epsilon$, $C_{\epsilon}$, is equal to
$Ch(R\Gamma_{U_{\epsilon}}{\cal G})$.
Moreover,
$$
\lim_{\epsilon \to 0^+} C_{\epsilon}=
\lim_{\epsilon \to 0^+} Ch(R\Gamma_{U_{\epsilon}}{\cal G})=
Ch(R\Gamma_U{\cal G}) = Ch \bigl( (Rj_{U \hookrightarrow M})_* \circ
(j_{U \hookrightarrow M})^* ({\cal G}) \bigr).
$$
\end{prop}

Proof of this proposition is completely analogous.

\separate

Recall that, in general, if $\tilde C_{(0,a)}$ is a family of cycles
in some space $Z$, then
$\lim_{\epsilon \to 0^+} \tilde C_{\epsilon} = \tilde C_0$
means that if we regard $\tilde C_{(0,a)}$ as a chain in $[0,a) \times Z$,
then $\partial \tilde C_{(0,a)} = - \tilde C_0$.
Thus it is natural to define $\lim_{\epsilon \to a^-} \tilde C_{\epsilon}$
by setting
$$
\lim_{\epsilon \to a^-} \tilde C_{\epsilon} =
\partial \tilde C_{(0,a)},
$$
where $\tilde C_{(0,a)}$ is regarded as a chain in $(0,a] \times Z$.

\separate

The following observation is completely trivial, but will play a very
important role.

\begin{rem} \label{deformationremark}
Let $C_{(0,a)}$ be a family of cycles in $T^*M$ as in Proposition
\ref{deformation1} or \ref{deformation2}.
Recall that a family of cycles $C_{(0,a)}$ is a cycle in
$(0,a) \times T^*M$.
We can regard it as a chain in $[0,a] \times T^*M$.
Let $p: [0,a] \times T^*M \twoheadrightarrow T^*M$ be the projection.
Then $p_*(C_{(0,a)})$ is a conic Borel-Moore
$(\dim_{\BB R}M+1)$-chain in $T^*M$ such that
$$
\partial (p_*(C_{(0,a)}))
= \lim_{\epsilon \to a^-} C_{\epsilon} - \lim_{\epsilon \to 0^+} C_{\epsilon}.
$$
In particular, the cycles $\lim_{\epsilon \to 0^+} C_{\epsilon}$
and $\lim_{\epsilon \to a^-} C_{\epsilon}$ are homologous.
In the situation of Proposition \ref{deformation1}, the support
of $p_*(C_{(0,a)})$ lies in the closure
$$
\overline{ \bigl\{ |Ch{\cal G}| \cap T^*U - sdf; \, s \ge 0 \bigr\}}.
$$
Similarly, in the situation of Proposition \ref{deformation2},
the support of $p_*(C_{(0,a)})$ lies in
$$
\overline{ \bigl\{ |Ch{\cal G}| \cap T^*U + sdf; \, s \ge 0 \bigr\}}
$$
(notice the difference in signs of $s df$).
\end{rem}

\begin{rem}
Let $C_{(0,a)}$ be a family of cycles in $T^*M$ as in Proposition
\ref{deformation1} (respectively Proposition \ref{deformation2}).
Suppose that there is $b > a$ such that the condition
$|Ch({\cal G})| \cap df \cap T^* \widetilde{V} = \varnothing$,
(respectively $|Ch({\cal G})| \cap -df \cap T^* \widetilde{V} = \varnothing$)
is satisfied on a bigger open set
$\widetilde{V} = \{ u \in U; 0 < f(u) < b \}$, then it is easy to show that
$$
\lim_{\epsilon \to a^-} C_{\epsilon} = C_a.
$$
But in general, $\lim_{\epsilon \to a^-} C_{\epsilon}$ need not equal $C_a$.
\end{rem}

This construction is similar to the classical Morse's lemma.
Recall that it says
that if we have a smooth real valued function $f$ on a manifold
$M$, then the sublevel sets $\{ m \in M;\, f(m) < a \}$ and
$\{ m \in M;\, f(m) < b \}$ can be deformed one into the other
as long as there are no critical values of $f$ in an open
interval containing $a$ and $b$.

\separate

\begin{section}
{Deformation of $Ch({\cal F})$ in $T^*X$}  \label{deformationsection}
\end{section}

Recall that we view ${\cal F}$ not as an element of the
``$G_{\BB R}$-equivariant derived category on $X$ with twist
$(-\lambda - \rho)$'' denoted in \cite{SchV2} by
$\operatorname{D}_{G_{\BB R}}(X)_{-\lambda}$,
but simply as a $G_{\BB R}$-equivariant sheaf on
the flag variety $X$ with the same characteristic cycle.
In this section we start with an element $g_0 \in \g g_{\BB R}'$
and the characteristic cycle $Ch({\cal F})$
of a $G_{\BB R}$-equivariant sheaf ${\cal F}$ on the flag variety $X$
and use general results of the previous section to deform
$Ch({\cal F})$ into a cycle of the form
$$
m_1 T^*_{x_1}X +\dots+ m_{|W|} T^*_{x_{|W|}}X,
$$
where $m_1,\dots,m_{|W|}$ are some integers,
$x_1,\dots,x_{|W|}$ are the points in $X$ fixed by $g_0$,
and each cotangent space $T^*_{x_k}X$ is given some orientation.
Moreover, to ensure good behavior of our integral (\ref{theintegral}),
we will stay during the process of deformation inside the set
\begin{equation}  \label{deformationset}
\{\zeta \in T^*X;\, Re( \langle g_0,\mu(\zeta) \rangle) \le 0\}.
\end{equation}
The precise result is stated in Proposition \ref{C(g)prop}.
This deformation will help us calculate the integral (\ref{theintegral}).

\separate

\begin{rem}  \label{triangle1}
Suppose ${\cal G}$ is a sheaf on a manifold $M$ and $Z$ is a
locally closed subset of $M$. Let $i: Z \hookrightarrow M$ be
the inclusion. Then M. Kashiwara and P. Schapira introduce
in \cite{KaScha}, Chapter II, a sheaf $i_! \circ i^* ({\cal G})$
denoted by ${\cal G}_Z$.
If $Z'$ is closed in $Z$, then they prove existence of a
distinguished triangle
$$
{\cal G}_{Z \setminus Z'} \to {\cal G}_Z \to {\cal G}_{Z'}.
$$
Hence
$$
Ch( {\cal G}_Z ) = Ch( {\cal G}_{Z \setminus Z'} )
+ Ch( {\cal G}_{Z'}).
$$
\end{rem}

Pick an element $g_0 \in \g g_{\BB R}'$ and let
$\g t \subset \g g$ be the unique (complex) Cartan subalgebra
containing $g_0$. Let $\Psi \subset \g t^*$ be the root system
of $\g g$ with respect to $\g t$.
Pick a positive root system $\Psi^{\le} \subset \Psi$ such that
$Re(\alpha(g_0)) \le 0$ for all $\alpha \in \Psi^{\le}$.
Let $B_{g_0} \subset G$ be the Borel subgroup whose Lie algebra
contains $\g t$ and the positive root spaces.
The action of $B_{g_0}$ on $X$ has $|W|$ orbits $O_1, \dots, O_{|W|}$.
Each orbit $O_k$ is a locally closed subset.
Hence it follows from Remark \ref{triangle1} that, as an element
of the $K$-group of the category of bounded complexes of
$\BB R$-constructible sheaves on $X$, our sheaf $\cal F$ is equivalent to
${\cal F}_{O_1}+ \dots + {\cal F}_{O_{|W|}}$,
and so
$$
Ch({\cal F})=Ch({\cal F}_{O_1})+ \dots + Ch({\cal F}_{O_{|W|}}).
$$

The idea is to deform each summand $Ch({\cal F}_{O_k})$
separately.
Let $x_k \in X$ be the only point in $O_k$ fixed by $g_0$;
it determines a Borel subalgebra $\g b_k \subset \g g$ 
consisting of all elements of $\g g$ fixing $x_k$.
This Borel subalgebra $\g b_k$ in turn determines
a positive root system (different from $\Psi^{\le}$)
such that $\g b_k$ contains $\g t$ and all the negative root spaces.
Let $\g n_k$ be the nilpotent subalgebra of $\g g$ containing all
the positive root spaces with respect to $\g b_k$, so that
$\g g = \g b_k \oplus \g n_k$ as linear spaces.
Define a subset $\Psi^{\le}_k \subset \Psi^{\le}$ by
$$
\Psi^{\le}_k = \{ \alpha \in \Psi^{\le};\, \g n_{\alpha} \subset \g n_k\},
$$
where $\g n_{\alpha}$ is the $\alpha$-root space in $\g g$.
We break $\Psi^{\le}_k$ into two subsets:
$$
\Psi^<_k = \{\alpha \in \Psi^{\le}_k; \, Re(\alpha(g_0)) < 0 \}
\quad \text{and} \quad
\Psi^0_k = \{\alpha \in \Psi^{\le}_k; \, Re(\alpha(g_0)) = 0 \}.
$$
Next we define three nilpotent subalgebras of $\g n_k$:
$$
\g n^<_k = \bigoplus_{\alpha \in \Psi^<_k} \g n_{\alpha},
\quad
\g n^0_k = \bigoplus_{\alpha \in \Psi^0_k} \g n_{\alpha}
\quad \text{and} \quad
\g n^>_k =
\bigoplus_{\alpha \notin \Psi^{\le}_k, \, \g n_{\alpha} \subset \g n_k}
\g n_{\alpha}.
$$
Of course, $\g n_k = \g n^<_k \oplus \g n^0_k \oplus \g n^>_k$
as linear spaces.

We define a map
$\psi_{g_0,k}: \g n_k \to X$,
$n \mapsto \exp(n) \cdot x_k$.
Then $\psi_{g_0,k}$ is a diffeomorphism of $\g n_k$ onto its image.
Let $U_k = \psi_{g_0,k} (\g n_k)$ be this image.

We decompose the nilpotent algebra $\g n_k$ into root spaces with respect to $\g t$:
$\g n_k= \bigoplus_{l=1}^n \g n_{\alpha_{x_k,l}}$.
We will assume that the roots $\alpha_{x_k,1},\dots,\alpha_{x_k,n}$
are enumerated so that
$\g n^0_k= \bigoplus_{l=1}^m \g n_{\alpha_{x_k,l}}$.
Each $\g n_{\alpha_{x_k,l}}$ is a one-dimensional complex vector space,
so we can
choose a linear coordinate $z_l: \g n_{\alpha_{x_k,l}} \tilde \to \BB C$.
Define a norm $\|.\|_k$ on $\g n_k$ by
$\|z\|_k = \sqrt{|z_1|^2+\dots+|z_n|^2}$.

Recall that we denote by $\operatorname{VF}_g$ the vector field on $X$
generated by $g \in \g g$: if $x \in X$ and $f \in {\cal C}^{\infty}(X)$, then
$$
\operatorname{VF}_g (x)f=
\frac d{d\epsilon} f(\exp(\epsilon g) \cdot x)|_{\epsilon=0}.
$$
Then
\begin{equation}  \label{mu}
\langle g, \mu(\zeta) \rangle  = \langle \operatorname{VF}_g, \zeta \rangle.
\end{equation}

\begin{lem} \label{linearvf}
For each $\tilde g \in \g t$, the vector field $\operatorname{VF}_{\tilde g}$
is expressed in the coordinate system $\psi_{g_0,k}$ by
\begin{equation}   \label{vfield}
\alpha_{x_k,1}(\tilde g) z_1 \frac{\partial}{\partial z_1} +\dots+
\alpha_{x_k,n}(\tilde g) z_n \frac{\partial}{\partial z_n}.
\end{equation}
\end{lem}

\pf
Let $T \subset G $ be the connected Lie subgroup
generated by $\g t$; $T$ is a Cartan subgroup of $G$.
First we observe that the map $\psi_{g_0,k} : \g n_k \to X$ is
$T$-equivariant. Indeed, if $n \in \g n_k$ and $t \in T$, then
\begin{multline*}
\psi_{g_0,k} (t \cdot n) = \exp \bigl( Ad(t)(n) \bigr) \cdot x_k  \\
= (t \exp(n) t^{-1}) \cdot x_k = (t \exp(n)) \cdot x_k
= t \cdot \psi_{g_0,k}(n).
\end{multline*}

For $n \in \g n_k = \bigoplus_{l=1}^n \g n_{\alpha_{x_k,l}}$,
say $n=(n_{\alpha_{x_k,1}}, \dots, n_{\alpha_{x_k,n}})$,
we obtain using the $T$-equivariance of $\psi_{g_0,k}$ that
\begin{multline*}
e^{\epsilon \tilde g} \cdot \psi_{g_0,k} (n)
= \psi_{g_0,k} \bigl( Ad(e^{\epsilon \tilde g})(n_{\alpha_{x_k,1}}, \dots,
n_{\alpha_{x_k,n}}) \bigr)  \\
= \psi_{g_0,k} (e^{\alpha_{x_k,1}(\epsilon \tilde g)} n_{\alpha_{x_k,1}},
\dots, e^{\alpha_{x_k,n}(\epsilon \tilde g)} n_{\alpha_{x_k,n}}).
\end{multline*}
That is,
$$
e^{\epsilon \tilde g} \cdot \psi_{g_0,k} (z_1,\dots,z_n) =
\psi_{g_0,k} (e^{\epsilon\alpha_{x_k,1}(\tilde g)} z_1, \dots, 
e^{\epsilon\alpha_{x_k,n}(\tilde g)} z_n);
$$
and the lemma follows from the definition of $\operatorname{VF}_{\tilde g}$.
\qed

Observe that the restriction $\psi_{g_0,k}|_{\g n^<_k \oplus \g n^0_k}$
is a diffeomorphism of $\g n^<_k \oplus \g n^0_k$ onto $O_k$.
Let $N^0_k = \psi_{g_0,k}(\g n^0_k)$, it is a closed subset of $O_k$.

If $g_0$ lies in a split Cartan $\g t_{\BB R} \subset \g g_{\BB R}$,
then $\g n^0_k=0$ and $N^0_k=\{x_k\}$.
On the other extreme, if the Cartan subalgebra
$\g t_{\BB R} \subset \g g_{\BB R}$ containing $g_0$ is compact,
then $\g n^0_k=\g n_k$ and $N^0_k=O_k$.

\separate

\begin{lem}  \label{open}
The $G_{\BB R}$-orbit of $x_k$ intersected with $N^0_k$
(i.e. $G_{\BB R} \cdot x_k \cap N^0_k$) contains an open neighborhood
of $x_k$ in $N^0_k$.
\end{lem}

\pf
Recall that $\g b_k$ denotes the Borel subalgebra in $\g g$ consisting
of all elements of $\g g$ fixing $x_k$ and that we have
$\g g = \g b_k \oplus \g n_k$ as linear spaces.
The tangent space of the flag variety
$X$ at $x_k$ can be naturally identified with $\g g / \g b_k \simeq \g n_k$.
It follows that to prove that $G_{\BB R} \cdot x_k \cap N^0_k$
contains an open neighborhood of $x_k$ in $N^0_k$ it is enough
to show that $\g g_{\BB R} + \g b_k \supset \g n^0_k$.

So pick a root $\alpha \in \Psi$ such that the root space
$\g n_{\alpha} \subset \g n^0_k$, i.e. $Re(\alpha(g_0))=0$.
By definition, $\g g_{\BB R}'$ consists of regular semisimple elements
$\tilde g \in \g g_{\BB R}$ which satisfy the following additional property:
If $\tilde{\g t}_{\BB R} \subset \g g_{\BB R}$ and
$\tilde{\g t} \subset \g g$ are the unique Cartan subalgebras in
$\g g_{\BB R}$ and $\g g$ respectively containing $\tilde g$,
and $\beta \in \tilde{\g t}^*$ is a (complex) root such that
$Re(\beta)|_{\tilde{\g t}_{\BB R}} \not \equiv 0$,
then $Re(\beta(\tilde g)) \ne 0$.
Hence $Re(\alpha(g_0))=0$ implies that
$Re(\alpha)|_{\g t_{\BB R}} \equiv 0$.

The Lie subalgebras
$\g n_{\alpha} \oplus \g n_{-\alpha} \oplus [\g n_{\alpha},\g n_{-\alpha}]$
and
$\g g_{\BB R} \cap
(\g n_{\alpha} \oplus \g n_{-\alpha} \oplus [\g n_{\alpha},\g n_{-\alpha}])$
of $\g g$ are isomorphic to $\g {sl}(2,\BB C)$ and $\g {sl}(2,\BB R)$
respectively.
Then $Re(\alpha)|_{\g t_{\BB R}} \equiv 0$ implies that
$\tilde{\g t}_{\BB R} = \g g_{\BB R} \cap
(\g n_{-\alpha} \oplus [\g n_{\alpha},\g n_{-\alpha}])$
is a compact Cartan subalgebra of
$\g g_{\BB R} \cap
(\g n_{\alpha} \oplus \g n_{-\alpha} \oplus [\g n_{\alpha},\g n_{-\alpha}])
\simeq \g {sl}(2,\BB R)$.
It is a well-known property of $\g {sl}(2,\BB R)$ that if
$\tilde{\g t}_{\BB R}$ is a compact Cartan of $\g {sl}(2,\BB R)$
and $\tilde {\g b} \subset \g {sl}(2,\BB C)$ is a complex Borel subalgebra
containing $\tilde{\g t}_{\BB R}$, then
$$
\g {sl}(2,\BB R) + \tilde {\g b} = \g {sl}(2,\BB C).
$$
In particular,
$$
\bigl( \g g_{\BB R} \cap
(\g n_{\alpha} \oplus \g n_{-\alpha} \oplus [\g n_{\alpha},\g n_{-\alpha}])
\bigr)
+ \bigl (\g n_{-\alpha} \oplus [\g n_{\alpha},\g n_{-\alpha}] \bigr)
\supset \g n_{\alpha},
$$
and $\g g_{\BB R} + \g b_k \supset \g n_{\alpha}$.
Since this is true for all root spaces $\g n_{\alpha} \subset \g n^0_k$,
$\g g_{\BB R} + \g b_k \supset \g n^0_k$.
\qed

This lemma implies that there exists an $r>0$ such that
\begin{equation}  \label{r}
\psi_{g_0,k} \bigl(\{ n^0 \in \g n^0_k;\, \|n^0\|_k < 2 r\} \bigr)
\subset (G_{\BB R} \cdot x_k \cap N^0_k).
\end{equation}
Let $B_k$ be the image under $\psi_{g_0,k}$ of the
open cylinder of radius $r$ in $\g n_k^< \oplus \g n_k^0$:
$$
B_k = \psi_{g_0,k} \bigl(
\{ (n^<,n^0,n^>) \in \g n^<_k \oplus \g n^0_k \oplus \g n^>_k = \g n_k;\,
n^> =0, \, \|n^0\|_k < r \} \bigr).
$$
$B_k$ is an open subset of $O_k$, and if $\g n_k^0 =0$
(i.e. the Cartan algebra $\g t_{\BB R} \subset \g g_{\BB R}$
containing $g_0$ is split), then $B_k = O_k$.

According to Remark \ref{triangle1} we have a distinguished triangle:
$$
{\cal F}_{B_k} \to {\cal F}_{O_k} \to {\cal F}_{O_k \setminus B_k},
$$
and hence
$$
Ch({\cal F}_{O_k}) = Ch({\cal F}_{B_k})
+ Ch({\cal F}_{O_k \setminus B_k}).
$$

\separate

Recall that $\g t_{\BB R} \subset \g g_{\BB R}$ is the Cartan subalgebra
containing $g_0$.
Let $T_{\BB R} \subset G_{\BB R}$ be the connected Lie subgroup
generated by $\g t_{\BB R}$. Notice that because we require $T_{\BB R}$
to be connected it may not be a Cartan subgroup of $G_{\BB R}$.
Recall that the sheaf ${\cal F}$ is $G_{\BB R}$-equivariant.
Hence $Ch({\cal F})$ is $T_{\BB R}$-invariant, and so
$$
Re( \langle g,\mu(\zeta) \rangle) =
Re( \langle \operatorname{VF}_g, \zeta \rangle) = 0
$$
for all $g \in \g t_{\BB R}$ and all $\zeta \in |Ch({\cal F})|$.

Similarly, because each of the sets $O_k$, $B_k$ and
$O_k \setminus B_k$ is $T_{\BB R}$-invariant, the sheaves
${\cal F}_{O_k}$, ${\cal F}_{B_k}$ and
${\cal F}_{O_k \setminus B_k}$ are $T_{\BB R}$-equivariant,
their characteristic cycles are $T_{\BB R}$-invariant, and
$Re( \langle g,\mu(\zeta) \rangle) =0$ for all 
$g \in \g t_{\BB R}$ and all $\zeta$ in
$|Ch( {\cal F}_{O_k} )|$, $|Ch( {\cal F}_{B_k} )|$ and
$|Ch( {\cal F}_{O_k \setminus B_k})|$.

\begin{lem}
The cycle $Ch({\cal F}_{O_k \setminus B_k})$
is homologous to the zero cycle inside the set
$\{\zeta \in T^*X;\, Re (\langle g_0, \mu(\zeta) \rangle) \le 0 \}$.
\end{lem}

\pf
If $\g n_k^0 =0$, then $O_k \setminus B_k = \varnothing$
and we are done. So let us suppose that $\g n_k^0 \ne 0$.

The sheaf ${\cal F}_{O_k \setminus B_k}$ is the extraordinary
direct image of a sheaf on $U_k$:
$$
{\cal F}_{O_k \setminus B_k} = (j_{U_k \hookrightarrow X})_! \circ
(j_{O_k \setminus B_k \hookrightarrow U_k})_! ({\cal F}|_{O_k \setminus B_k}).
$$

Let $f: U_k \to \BB R$ be the function defined by
$f(\exp(n) \cdot x_k) = e^{-\|n\|_k^2}$.
This function extends by zero to a smooth function on all of $X$.
It follows from Lemma \ref{linearvf} that, for each $x \in O_k$,
$Re(\langle g_0, \mu(df(x)) \rangle) =
Re( \langle \operatorname{VF}_{g_0}, df(x) \rangle) \ge 0$.

By the open embedding theorem,
\begin{multline*}
Ch( {\cal F}_{O_k \setminus B_k} ) =
Ch \bigl( (j_{U_k \hookrightarrow X})_! \circ 
(j_{O_k \setminus B_k \hookrightarrow U_k})_!
({\cal F}|_{O_k \setminus B_k}) \bigr)  \\
= \lim _{s \to 0^+}
Ch \bigl( (j_{O_k \setminus B_k \hookrightarrow U_k})_!
({\cal F}|_{O_k \setminus B_k}) \bigr) - s\frac {df}f.
\end{multline*}
In particular, there exists a chain $C$ in $T^*X$ such that
$$
\partial C =
Ch( {\cal F}_{O_k \setminus B_k} ) -
\Bigl( Ch \bigl( (j_{O_k \setminus B_k \hookrightarrow U_k})_!
({\cal F}|_{O_k \setminus B_k}) \bigr) - \frac {df}f \Bigr)
$$
and the support of this chain lies inside the set (\ref{deformationset}).
Notice that
$$
Ch \bigl( (j_{O_k \setminus B_k \hookrightarrow U_k})_!
({\cal F}|_{O_k \setminus B_k}) \bigr) - \frac {df}f
$$
is a cycle in $T^*X$ whose support lies inside $T^*U$.

Recall that the roots $\alpha_{x_k,1},\dots,\alpha_{x_k,n}$
are enumerated so that
$\g n^0_k= \bigoplus_{l=1}^m \g n_{\alpha_{x_k,l}}$, i.e.
the complex numbers
$\alpha_{x_k,1}(g_0),\dots,\alpha_{x_k,m}(g_0)$ are purely
imaginary (and nonzero).
So let us write
$$
\alpha_{x_k,1}(g_0)=i \tilde \alpha_{x_k,1}(g_0),\dots,
\alpha_{x_k,m}(g_0)=i \tilde \alpha_{x_k,m}(g_0).
$$
If we write each $z_l$ as $x_l + iy_l$, then the formula (\ref{vfield})
in Lemma \ref{linearvf} expressing the vector field
$\operatorname{VF}_{g_0}$ on $X$ generated by $g_0$
in the coordinate system $\psi_{g_0,k}$ becomes
\begin{multline*}
\tilde \alpha_{x_k,1}(g_0) \bigl( - y_1 \frac{\partial}{\partial x_1} +
x_1 \frac{\partial}{\partial y_1} \bigr) + \dots +
\tilde \alpha_{x_k,m}(g_0) \bigl( - y_m \frac{\partial}{\partial x_m} +
x_m \frac{\partial}{\partial y_m} \bigr)   \\
+ \beta_{m+1}(x_{m+1},y_{m+1}) \frac{\partial}{\partial x_{m+1}} +
\beta_{m+1}'(x_{m+1},y_{m+1}) \frac{\partial}{\partial y_{m+1}} + \dots  \\
+ \beta_n(x_n,y_n) \frac{\partial}{\partial x_n} +
\beta_n'(x_n,y_n) \frac{\partial}{\partial y_n},
\end{multline*}
for some real coefficients
$$
\beta_{m+1}(x_{m+1},y_{m+1}), \beta_{m+1}'(x_{m+1},y_{m+1}), \dots,
\beta_n(x_n,y_n),\beta_n'(x_n,y_n).
$$
which depend linearly on respective $(x_l,y_l)$ and also depend on $g_0$.
Define a 1-form $\eta$ on $\g n_k$ to be
$$
\eta  =
\frac {\tilde \alpha_1(g_0)}{|\tilde \alpha_1(g_0)|}
(y_1 dx_1 - x_1 dy_1) +\dots+
\frac {\tilde \alpha_m(g_0)}{|\tilde \alpha_m(g_0)|} (y_m dx_m - x_m dy_m).
$$
We can regard $\eta$ as a form on $U_k$ via the diffeomorphism $\psi_{g_0,k}$.
Then
$$
Re \bigl( \bigl\langle g_0, \mu(\eta(\psi_{g_0,k}(n^<,n^0,n^>)))
\bigr\rangle \bigr) =
- |\tilde \alpha_1(g_0)| (x_1^2 + y_1^2) - \dots
- |\tilde \alpha_m(g_0)| (x_m^2 + y_m^2)
$$
which is at most
$$
-r^2 \min\{|\tilde \alpha_1(g_0)|,\dots,|\tilde \alpha_m(g_0)|\}
$$
on $O_k \setminus B_k$.

Finally, define a $(2n+1)$-chain in $T^*X$
$$
\tilde C = - \Bigl(
Ch \bigl( (j_{O_k \setminus B_k \hookrightarrow U_k})_!
({\cal F}|_{O_k \setminus B_k}) \bigr) - \frac {df}f \Bigr)
+ t \eta, \qquad t \in [0,\infty).
$$
Then
$$
\partial \tilde C =
Ch \bigl( (j_{O_k \setminus B_k \hookrightarrow U_k})_!
({\cal F}|_{O_k \setminus B_k}) \bigr) - \frac {df}f
$$
and the support $|\tilde C|$ lies in the set (\ref{deformationset}).
\qed

\separate

Our next task is to deform $Ch({\cal F}_{B_k})$.
Let ${\cal G}={\cal F}_{B_k} =
(j_{B_k \hookrightarrow X})_! \circ (j_{B_k \hookrightarrow X})^* ({\cal F})$.

Let $\gamma: [0,\infty) \to [1,\infty)$ be a constructible
${\cal C}^2$-function such that $\gamma([0,4r^2]) = \{1\}$ and
$\gamma(t) = t$ for $t$ sufficiently large.
We define a function $f: U_k \to \BB R$ by
\begin{equation}  \label{f}
f(\psi_{g_0,k}(n^<+n^0+n^>))
= \exp \bigl( -(\gamma(\|n^0\|_k^2) - \|n^<\|_k^2 - \|n^>\|_k^2 \bigr)
\end{equation}
and extend it by zero to a smooth function on all of $X$.
It follows from Lemma \ref{linearvf} that, for each $x$ in the closure
$\overline{B_k}$,
$Re(\langle g_0, \mu(df(x)) \rangle) =
Re( \langle \operatorname{VF}_{g_0}, df(x) \rangle) \ge 0$.
Moreover, if we define
$$
V = \{ u \in U_k;\, 0<f(u)<1 \} = U_k \setminus f^{-1}(1),
$$
then $Re(\langle \operatorname{VF}_{g_0}, dh \rangle)$ is strictly
positive on $V \cap \overline{B_k}$.
Since $\supp({\cal G}) \subset \overline{B_k}$, the condition
$|Ch({\cal G})| \cap df \cap T^*V = \varnothing$
of Proposition \ref{deformation1} with $a=1$ is satisfied.
Thus we can apply Proposition \ref{deformation1} with ambient
manifold $X$, open subset $U_k = \psi_{g_0,k}(\g n_k)$, function
$f$ defined by (\ref{f})
to get a family of cycles $C_{(0,1)}$ such that,
for each $\epsilon \in (0,1)$,
$$
C_{\epsilon}=Ch({\cal G}_{U_{\epsilon}}) \qquad \text{and} \qquad
\lim_{\epsilon \to 0^+} C_{\epsilon} = Ch({\cal G}_{U_k}) = Ch({\cal G}),
$$
where $U_{\epsilon}=\{u \in U_k; \, f(u)>\epsilon\}$.

Remark \ref{deformationremark} tells us that the difference of cycles
$Ch({\cal G})$ and $\lim_{\epsilon \to 1^-} C_{\epsilon}$ is the
boundary of a certain chain. Because
$Re(\langle g_0, \mu(df) \rangle) >0$ on $V \cap \overline{B_k}$,
this chain is supported inside the set (\ref{deformationset}).

\separate

We need to determine $\lim_{\epsilon \to 1^-} C_{\epsilon}$.
Let $\tilde C_{(0,1)}$ be the family of cycles $C_{(0,1)}$ restricted
to $T^*U_k$. Then, for each $\epsilon \in (0,1)$,
$$
\tilde C_{\epsilon} =
Ch \bigl( (j_{U_\epsilon \hookrightarrow U_k})_!
({\cal G}|_{U_{\epsilon}}) \bigr)
\qquad \text{and} \qquad
C_{\epsilon} = Ch({\cal G}_{U_{\epsilon}}) =
(j_{U_k \hookrightarrow X})_! \tilde C_{\epsilon}.
$$
Notice that when $\epsilon > 1/2$ the support of $C_{\epsilon}$ lies
inside $\pi^{-1}(\{x \in X;\, f(x) \ge 1/2\})$, and the set
$\{x \in X;\, f(x) \ge 1/2\}$ is a compact subset of $U_k$.
This implies that
\begin{multline*}
\lim_{\epsilon \to 1^-} Ch({\cal G}_{U_{\epsilon}})=
\lim_{\epsilon \to 1^-} C_{\epsilon}=
(j_{U_k \hookrightarrow X})_!(\lim_{\epsilon \to 1^-} \tilde C_{\epsilon})  \\
= (j_{U_k \hookrightarrow X})_! \bigl( \lim_{\epsilon \to 1^-}
Ch \bigl( (j_{U_{\epsilon} \hookrightarrow U_k})_!({\cal G}|_{U_{\epsilon}})
\bigr) \bigr).
\end{multline*}

Thus we need to find
$$
\lim_{\epsilon \to 1^-} \tilde C_{\epsilon}=
\lim_{\epsilon \to 1^-} Ch \bigl(
(j_{U_{\epsilon} \hookrightarrow U_k})_!({\cal G}|_{U_{\epsilon}}) \bigr).
$$

Let $V_{\epsilon}= \{ u \in U_k;\, 0 < f(u) < \epsilon \}$, so that
$V_1=V=\{ u \in U_k;\, 0<f(u)<1 \}$.

\separate

\begin{claim}
$$
Ch \bigl(
(j_{U_{\epsilon} \hookrightarrow U_k})_!({\cal G}|_{U_{\epsilon}}) \bigr) =
Ch({\cal G}|_{U_k}) -
Ch \bigl( (R j_{V_{\epsilon} \hookrightarrow U_k})_*
({\cal G}|_{V_{\epsilon}}) \bigr).
$$
\end{claim}

\pf
Notice that 
$$
\bigl( (j_{U_{\epsilon} \hookrightarrow U_k})_!
({\cal G}|_{U_{\epsilon}}) \bigr) |_{U_{\epsilon}}
= ({\cal G}|_{U_k}) |_{U_{\epsilon}}
\qquad \text{and} \qquad
\bigl( (R j_{V_{\epsilon} \hookrightarrow U_k})_*
({\cal G}|_{V_{\epsilon}}) \bigr) |_{U_{\epsilon}}= 0.
$$
Similarly,
$$
\bigl( (j_{U_{\epsilon} \hookrightarrow U_k})_!
({\cal G}|_{U_{\epsilon}}) \bigr) |_{V_{\epsilon}} =0
\qquad \text{and} \qquad
({\cal G}|_{U_k}) |_{V_{\epsilon}} =
\bigl( (R j_{V_{\epsilon} \hookrightarrow U_k})_*
({\cal G}|_{V_{\epsilon}}) \bigr) |_{V_{\epsilon}}.
$$
This shows that the difference of these two cycles,
$$
Ch({\cal G}|_{U_k}) -
Ch \bigl( (R j_{V_{\epsilon} \hookrightarrow U_k})_*
({\cal G}|_{V_{\epsilon}}) \bigr) -
Ch \bigl( (j_{U_{\epsilon} \hookrightarrow U_k})_!
({\cal G}|_{U_{\epsilon}}) \bigr),
$$
is supported in
\begin{multline*}
\pi^{-1} (\supp({\cal G})) \cap
\bigl( T^*U_k \setminus T^*(U_{\epsilon} \cup V_{\epsilon}) \bigr)  \\
= \pi^{-1} \bigl( f^{-1}(\epsilon) \cap \supp({\cal G}) \bigr) \subset
\pi^{-1} (f^{-1}(\epsilon) \cap \overline{B_k}).
\end{multline*}
Thus it is enough to show that each
$x \in f^{-1}(\epsilon) \cap \overline{B_k}$
lies in an open neighborhood $\Omega_x \subset U_k$
such that the cycle in $T^*\Omega_x$
\begin{multline}  \label{zero}
Ch \bigl( ({\cal G}|_{U_k}) |_{\Omega_x} \bigr) -
Ch \bigl( \bigl( (R j_{V_{\epsilon} \hookrightarrow U_k})_*
({\cal G}|_{V_{\epsilon}}) \bigr) \bigl|_{\Omega_x} \bigr) -
Ch \bigl( \bigl( (j_{U_{\epsilon} \hookrightarrow U_k})_!
({\cal G}|_{U_{\epsilon}}) \bigr) \bigl|_{\Omega_x} \bigr)   \\
= Ch ({\cal G}|_{\Omega_x}) -
Ch \bigl( (R j_{V_{\epsilon} \cap \Omega_x \hookrightarrow \Omega_x})_*
({\cal G}|_{V_{\epsilon} \cap \Omega_x}) \bigr) -
Ch \bigl( (j_{U_{\epsilon} \cap \Omega_x \hookrightarrow \Omega_x})_!
({\cal G}|_{U_{\epsilon} \cap \Omega_x}) \bigr)
\end{multline}
is zero.

Because $Re(\langle \operatorname{VF}_{g_0}, df(x) \rangle) \ne 0$,
there exists an open neighborhood $\Omega_x \subset U_k$ containing $x$ and
smooth real-valued functions $y_2, \dots, y_{2n}$ defined on $\Omega_x$
such that $\{ f - \epsilon ,y_2, \dots, y_{2n} \}$ form a
${\cal C}^{\infty}$ system of coordinates on $\Omega_x$ centered
at $x$ and the vector field $\operatorname{VF}_{g_0}$ is expressed
in this coordinate system by
$$
\Bigl( F(f-\epsilon) \frac{\partial}{\partial(f-\epsilon)},0,\dots,0 \Bigr),
$$
for some nonvanishing single variable function $F$.
By making $\Omega_x$ smaller if necessary, we can assume that
$(f - \epsilon ,y_2, \dots, y_{2n})$ maps $\Omega_x$ diffeomorphically
onto an open cube $(-\delta,\delta)^{2n}$, for some $\delta >0$.

Until the rest of this proof we will regard
${\cal G}|_{\Omega_x}$ as a sheaf on $(-\delta,\delta)^{2n}$ via
this diffeomorphism.
Write $(-\delta,\delta)^{2n} = (-\delta,\delta) \times \widetilde{\Omega}_x$,
where $\widetilde{\Omega}_x = \{0\} \times (-\delta,\delta)^{2n-1}$,
and let $p_2: (-\delta,\delta)^{2n} \twoheadrightarrow \widetilde{\Omega}_x$
be the projection setting the first coordinate equal zero.
Then $\widetilde{\Omega}_x = f^{-1}(\epsilon) \cap \Omega_x$,
$V_{\epsilon} \cap \Omega_x = (-\delta,0) \times \widetilde{\Omega}_x$ and
$U_{\epsilon} \cap \Omega_x = (0,\delta) \times \widetilde{\Omega}_x$.
Also let
$\BB C_{(-\delta,\delta)}$, $\BB C_{(0,\delta)}$ and $\BB C_{(-\delta,0)}$
denote the constant sheaves on
$(-\delta,\delta)$, $(0,\delta)$ and $(-\delta,0)$ respectively.

Since the sheaf ${\cal G}|_{U_k}$ is $T_{\BB R}$-equivariant,
\begin{eqnarray*}
Ch({\cal G}|_{\Omega_x}) &=&
Ch \bigl( (p_2)^* ({\cal G}|_{\widetilde{\Omega}_x}) \bigr) =
Ch \bigl( \BB C_{(-\delta,\delta)} \boxtimes
{\cal G}|_{\widetilde{\Omega}_x} \bigr),  \\
Ch \bigl( (R j_{V_{\epsilon} \cap \Omega_x \hookrightarrow \Omega_x})_*
({\cal G}|_{V_{\epsilon} \cap \Omega_x}) \bigr) &=&
Ch \bigl( (Rj_{(-\delta,0) \hookrightarrow (-\delta,\delta)})_*
\BB C_{(-\delta,0)} \boxtimes
{\cal G}|_{\widetilde{\Omega}_x} \bigr),  \\
Ch \bigl( (j_{U_{\epsilon} \cap \Omega_x \hookrightarrow \Omega_x})_!
({\cal G}|_{U_{\epsilon} \cap \Omega_x}) \bigr) &=&
Ch \bigl( (j_{(0,\delta) \hookrightarrow (-\delta,\delta)})_!
\BB C_{(0,\delta)} \boxtimes
{\cal G}|_{\widetilde{\Omega}_x} \bigr).
\end{eqnarray*}
It is easy to see that
$$
Ch (\BB C_{(-\delta,\delta)}) -
Ch \bigl( (Rj_{(-\delta,0) \hookrightarrow (-\delta,\delta)})_*
\BB C_{(-\delta,0)} \bigr) -
Ch \bigl( (j_{(0,\delta) \hookrightarrow (-\delta,\delta)})_!
\BB C_{(0,\delta)} \bigr) = 0.
$$
Hence it follows that the cycle (\ref{zero}) is zero.
This proves the claim.
\qed

Thus
$$
\lim_{\epsilon \to 1^-} Ch \bigl(
(j_{U_{\epsilon} \hookrightarrow U_k})_!({\cal G}|_{U_{\epsilon}}) \bigr)
= Ch({\cal G}|_{U_k}) -
\lim_{\epsilon \to 1^-} Ch \bigl( (R j_{V_{\epsilon} \hookrightarrow U_k})_*
({\cal G}|_{V_{\epsilon}}) \bigr).
$$
Let us consider a function $f'=1-f$ on $U_k$. Then the set of $u \in U_k$
where $f'(u)$ is positive is precisely the set $V=\{ u \in U_k;\,0<f(u)<1 \}$.
Then Proposition \ref{deformation2} applied to manifold $U_k$,
open set $V$, function $f'$ and sheaf ${\cal G}|_{U_k}$
tells us that there is a family of cycles $\tilde C'_{(0,1)}$ in $T^*U_k$
such that, for each $\epsilon \in (0,1)$,
the specialization at $1-\epsilon$,
$$
\tilde C'_{1-\epsilon}=
Ch \bigl( (R j_{V_{\epsilon} \hookrightarrow U_k})_*
({\cal G}|_{V_{\epsilon}}) \bigr).
$$
Therefore,
$$
\lim_{\epsilon \to 1^-}
Ch \bigl( (R j_{V_{\epsilon} \hookrightarrow U_k})_*
({\cal G}|_{V_{\epsilon}}) \bigr) =
\lim_{\epsilon \to 0^+} \tilde C'_{\epsilon} =
Ch \bigl( (R j_{V_1 \hookrightarrow U_k})_*({\cal G}|_{V_1}) \bigr).
$$
Hence (recall that $V_1=V$)
$$
\lim_{\epsilon \to 1^-}
Ch \bigl(
(j_{U_{\epsilon} \hookrightarrow U_k})_!({\cal G}|_{U_{\epsilon}}) \bigr) =
Ch({\cal G}|_{U_k}) - Ch \bigl( (R j_{V \hookrightarrow U_k})_*
({\cal G}|_{V}) \bigr).
$$

\begin{rem}
If ${\cal G}$ is a constructible sheaf on a manifold $M$,
$Z$ a closed subset of $M$, $U=M \setminus Z$ its complement
and $i: Z \to M$, $j:U \to M$ the inclusion maps,
then we have a distinguished triangle
$$
(Ri)_* \circ i^! ({\cal G}) \to {\cal G} \to (Rj)_* \circ j^* ({\cal G}).
$$
Hence
$$
Ch( {\cal G} ) = Ch \bigl( (Ri)_* \circ i^! ({\cal G}) \bigr)
+ Ch \bigl( (Rj)_* \circ j^* ({\cal G}) \bigr).
$$
\end{rem}

We apply this remark to the ambient manifold $U_k$,
closed subset $Z= U_k \setminus V = f^{-1}(1)$, open subset $V$ and
sheaf ${\cal G}|_{U_k}$ to get
$$
Ch({\cal G}|_{U_k})=
Ch \bigl( (R j_{f^{-1}(1) \hookrightarrow U_k})_* \circ
(j_{f^{-1}(1) \hookrightarrow U_k})^! ({\cal G}|_{U_k}) \bigr)
+ Ch \bigl( (R j_{V \hookrightarrow U_k})_*({\cal G}|_{V}) \bigr).
$$
Hence
$$
\lim_{\epsilon \to 1^-}
Ch \bigl( (R j_{V_{\epsilon} \hookrightarrow U_k})_*
({\cal G}|_{V_{\epsilon}}) \bigr) =
Ch \bigl( (R j_{f^{-1}(1) \hookrightarrow U_k})_* \circ
(j_{f^{-1}(1) \hookrightarrow U_k})^! ({\cal G}|_{U_k}) \bigr).
$$

Let
$$
V_k = \psi_{g_0,k} \bigl(
\{ (n^<,n^0,n^>) \in \g n^<_k \oplus \g n^0_k \oplus \g n^>_k = \g n_k;\,
\|n^0\|_k < r \} \bigr),
$$
$V_k$ is an open subset of $X$, $V_k \cap O_k = B_k$, and $B_k$ is a closed
subset of $V_k$.

Letting $B_k^0 = f^{-1}(1) \cap B_k = N_k^0 \cap B_k$ we get
\begin{multline*}
\lim_{\epsilon \to 1^-} Ch({\cal G}_{U_{\epsilon}})=
(j_{U_k \hookrightarrow X})_! \bigl( \lim_{\epsilon \to 1^-}
Ch((R j_{V_{\epsilon} \hookrightarrow U_k})_*
({\cal G}|_{V_{\epsilon}})) \bigr)  \\
= Ch \bigl( (j_{U_k \hookrightarrow X})_! \circ
(R j_{f^{-1}(1) \hookrightarrow U_k})_* \circ
(j_{f^{-1}(1) \hookrightarrow U_k})^! ({\cal G}|_{U_k}) \bigr)  \\
= Ch \bigl( (j_{U_k \hookrightarrow X})_! \circ
(j_{f^{-1}(1) \hookrightarrow U_k})_! \circ
(j_{f^{-1}(1) \hookrightarrow V_k})^! \circ
(j_{V_k \hookrightarrow U_k})^! ({\cal G}|_{U_k}) \bigr)  \\
= Ch \bigl( (j_{f^{-1}(1) \hookrightarrow X})_! \circ
(j_{f^{-1}(1) \hookrightarrow V_k})^! ({\cal F}_{B_k}|_{V_k}) \bigr)  \\
= Ch \bigl( (Rj_{f^{-1}(1) \hookrightarrow X})_* \circ
(j_{f^{-1}(1) \hookrightarrow V_k})^! \circ
(R j_{B_k \hookrightarrow V_k})_* ({\cal F}|_{B_k}) \bigr)  \\
= Ch \bigl( (Rj_{f^{-1}(1) \hookrightarrow X})_* \circ
(Rj_{B^0_k \hookrightarrow f^{-1}(1)})_* \circ
(j_{B^0_k \hookrightarrow B_k})^! ({\cal F}|_{B_k}) \bigr)  \\
= Ch \bigl( (Rj_{B^0_k \hookrightarrow X})_* \circ
(j_{B^0_k \hookrightarrow B_k})^! ({\cal F}|_{B_k}) \bigr).
\end{multline*}

That is we have deformed the cycle $Ch({\cal F}_{O_k})$ into
$$
Ch \bigl( (Rj_{B^0_k \hookrightarrow X})_* \circ
(j_{B^0_k \hookrightarrow B_k})^! ({\cal F}|_{B_k}) \bigr).
$$

\separate

In the case when $\g t_{\BB R}$ is split, $B_k = O_k$, $N^0_k=\{x_k\}$
and $B^0_k =\{x_k\}$, so
\begin{multline*}
Ch \bigl( (Rj_{B^0_k \hookrightarrow X})_* \circ
(j_{B^0_k \hookrightarrow B_k})^! ({\cal F}|_{B_k}) \bigr) \\
= Ch \bigl( (R j_{\{x_k\} \hookrightarrow X})_* \circ
(j_{\{x_k\} \hookrightarrow O_k})^! ({\cal F|}_{O_k}) \bigr) =
m_k T^*_{x_k}X,
\end{multline*}
where
$$
m_k= \chi \bigl( R\Gamma_{\{x_k\}} ({\cal F}|_{O_k})_{x_k} \bigr),
$$
and the orientation of $T^*_{x_k}X$ is chosen so that
if we write each $z_l$ as $x_l + iy_l$,
$l=1,\dots,n$, then the $\BB R$-basis of $T^*_{x_k}X$
\begin{equation}  \label{orientation'}
\{ (\psi_{g_0,k})_*|_0(dx_1),(\psi_{g_0,k})_*|_0(dy_1),\dots,
(\psi_{g_0,k})_*|_0(dx_n),(\psi_{g_0,k})_*|_0(dy_n) \}
\end{equation}
is positively oriented.
Thus we have deformed $Ch({\cal F}_{O_k})$ into a cycle of desired type.

\separate

So let us assume that $\g t_{\BB R}$ is not split.
Because the sheaf ${\cal F}$ is $G_{\BB R}$-equivariant, it follows
from (\ref{r}) that the cycle
$Ch \bigl( (Rj_{B^0_k \hookrightarrow X})_* \circ
(j_{B^0_k \hookrightarrow B_k})^! ({\cal F}|_{B_k}) \bigr)$ is just
$m_k Ch \bigl( (Rj_{B_k^0 \hookrightarrow U_k})_* (\BB C_{B_k^0}) \bigr)$,
where $\BB C_{B_k^0}$ denotes the constant sheaf on $B_k^0$ and
\begin{equation}  \label{m}
m_k = \chi \bigl( (Rj_{B^0_k \hookrightarrow X})_* \circ
(j_{B^0_k \hookrightarrow B_k})^! ({\cal F}|_{B_k})_{x_k} \bigr).
\end{equation}

The set $B_k^0$ is $T_{\BB R}$-invariant which implies that
$Re (\langle g_0, \mu(\zeta) \rangle) =0$ for all
$\zeta \in \bigl |Ch \bigl(
(Rj_{B_k^0 \hookrightarrow U_k})_* (\BB C_{B_k^0}) \bigr) \bigr|$.
We define a function $f$ on $U_k = \psi_{g_0,k} (\g n_k)$ by
$$
f(\psi_{g_0,k}(n)) = r^2 - \|n\|_k^2.
$$
Notice that $Re (\langle g_0, \mu(df(x)) \rangle) =0$ for $x \in B_k^0$.
Then using Proposition \ref{deformation2} and Remark \ref{deformationremark}
we can deform the cycle
$Ch \bigl( (Rj_{B_k^0 \hookrightarrow U_k})_* (\BB C_{B_k^0}) \bigr)$
into the cycle
$Ch \bigl( (Rj_{\{x_k\} \hookrightarrow U_k})_* (\BB C_{\{x_k\}}) \bigr)
= T^*_{x_k} X$,
so that in the process of deformation we always stay inside the set
(\ref{deformationset}).
Here the orientation of $T^*_{x_k} X$ is chosen so that
the $\BB R$-basis (\ref{orientation'}) of $T^*_{x_k} X$
is positively oriented.

\separate

Thus we have deformed the cycle  $Ch ({\cal F}_{O_k})$
into $m_k T^*_{x_k}X$, and so we obtain a deformation of the cycle
$Ch({\cal F})=Ch({\cal F}_{O_1})+ \dots + Ch({\cal F}_{O_{|W|}})$
into $m_1 T^*_{x_1}X + \dots + m_{|W|} T^*_{x_{|W|}}X$.
The coefficient formula (\ref{m})
holds for all $g_0 \in \g g_{\BB R}'$, no matter whether
the Cartan subalgebra $\g t_{\BB R}$ containing $g_0$ is split or not.

We will show in Appendix \ref{appendix} that these coefficients
$m_k$'s coincide with coefficients $d_{g_0,x_k}$'s in \cite{SchV2}.

\separate

Let us summarize the result of our deformation.

\begin{prop}  \label{C(g)prop}
Pick an element $g_0 \in \g g_{\BB R}'$. Then there is a
Borel-Moore chain $C(g_0)$ in $T^*X$
of dimension $2n+1$ with the following properties:

(i) $C(g_0)$ is conic, i.e. invariant under the action of the multiplicative
group of positive reals on $T^*X$;

(ii) The support of $C(g_0)$ lies in the set
$\{\zeta \in T^*X; Re( \langle g_0, \mu(\zeta) \rangle ) \le 0 \}$;

(iii) Let $x_1, \dots, x_{|W|}$ be the fixed points of $g_0$ in $X$, then
there are integers $m_1, \dots, m_{|W|}$ such that
$$
\partial C(g_0) =
Ch({\cal F}) -(m_1 T^*_{x_1}X + \dots + m_{|W|} T^*_{x_{|W|}}X).
$$
More specifically,
$$
m_k = \chi \bigl( (Rj_{B^0_k \hookrightarrow X})_* \circ
(j_{B^0_k \hookrightarrow B_k})^! ({\cal F}|_{B_k})_{x_k} \bigr),
$$
and the orientation of $T^*_{x_k}X$ is chosen so that
if we write each $z_l$ as $x_l + iy_l$,
then the $\BB R$-basis of $T^*_{x_k}X$
\begin{equation}   \label{orientation}
\{ (\psi_{g_0,k})_*|_0(dx_1),(\psi_{g_0,k})_*|_0(dy_1),\dots,
(\psi_{g_0,k})_*|_0(dx_n),(\psi_{g_0,k})_*|_0(dy_n) \}
\end{equation}
is positively oriented.
\end{prop}

\begin{rem}  \label{C(g)remark}
Let $g \in \g t_{\BB R} \cap \g g_{\BB R}'$ be such that, for each root
$\alpha \in \Psi$, $Re(\alpha(g))$ has the same sign as $Re(\alpha(g_0))$
and if $Re(\alpha(g_0)) = 0$, then $Im(\alpha(g))$ has the same sign as
$Im(\alpha(g_0))$. Then we can choose the Borel subgroup $B_g \subset G$
equal $B_{g_0}$. In this case the chain $C(g)$ is identical to $C(g_0)$.
\end{rem}

\separate

\begin{section}
{Integration}
\end{section}

In this section we will compute the integral (\ref{theintegral})
first under the assumption that the form $\phi$ is compactly
supported in $\g g_{\BB R}'$ and then in general.

\separate

Pick an element $g_0$ lying in the support of $\phi$ and
let $\g t_{\BB R}$ be the Cartan subalgebra containing $g_0$.
There exists an open neighborhood $\Omega$ of $g_0$ in
$\g g_{\BB R}'$ and a smooth map $\omega: \Omega \to G_{\BB R}$
with the following two properties:

(i) $\omega|_{\Omega \cap \g t_{\BB R}} \equiv e$, the identity element of
$G_{\BB R}$;

(ii) For every $g \in \Omega$,
the conjugate Cartan subalgebra
$\omega(g) \g t_{\BB R} \omega(g)^{-1}$
contains $g$.

Making $\Omega$ smaller if necessary, we can assume that both
$\Omega$ and $\Omega \cap \g t_{\BB R}$ are connected.
Let $\g t = \g t_{\BB R} \oplus i \g t_{\BB R} \subset \g g$
be the complex Cartan subalgebra containing $g_0$.

\begin{rem}
One cannot deal with the integral (\ref{theintegral})
``one Cartan algebra at a time'' and avoid introducing
a map like $\omega$ because the limit
$$
\lim_{R \to \infty}
\int_{\g t_{\BB R} \times (Ch({\cal F}) \cap \{\|\zeta\|<R\})}
e^{\langle g, \mu_{\lambda}(\zeta) \rangle}
\phi(g) (-\sigma+\pi^* \tau_{\lambda})^n
$$
may not exist.
\end{rem}

\separate

From now on we will assume that the support of $\phi$ lies in
$\Omega$. The case $\supp(\phi) \subset \g g_{\BB R}'$
can be reduced to this special case
by a partition of unity argument.

\separate

Our biggest obstacle to making any deformation argument in order to
compute the integral (\ref{theintegral}) is that the integration
takes place over a cycle which is not compactly supported and Stokes'
theorem no longer applies.
In order to overcome this obstacle, we will construct a deformation
$\Theta_t: \Omega \times T^*X \to \Omega \times T^*X$, $t \in [0,1]$,
such that $\Theta_0$ is the identity map;
$$
Re \bigl( (\Theta_t)^*\langle g, \mu(\zeta) \rangle \bigr)
<  Re(\langle g, \mu(\zeta) \rangle)
$$
for $t>0$, $g \in \Omega$ and $\zeta \in T^*X$
which does not lie in the zero section (Lemma \ref{realpartlemma});
$\Theta_t$ essentially commutes with scaling the fiber of $T^*X$
(Lemma \ref{scaling}).
The last two properties will imply that the integral
$$
\int_{\g g_{\BB R} \times (Ch({\cal F}) \cap \{\|\zeta\|<R\})}
(\Theta_t)^* \bigl( e^{\langle g, \mu_{\lambda}(\zeta) \rangle}
\phi (-\sigma+\pi^* \tau_{\lambda})^n \bigr)
$$
converges absolutely for $t \in (0,1]$.
Finally, the most important property of $\Theta_t$ is stated
in Lemma \ref{slanting} which essentially says that we can
replace our integrand
$$
e^{\langle g, \mu_{\lambda}(\zeta) \rangle}
\phi (-\sigma+\pi^* \tau_{\lambda})^n
$$
with the pullback
$$
(\Theta_t)^* \bigl( e^{\langle g, \mu_{\lambda}(\zeta) \rangle}
\phi (-\sigma+\pi^* \tau_{\lambda})^n \bigr).
$$

\separate

We introduce $|W|$ coordinate systems on $\Omega \times T^*X$
in which the integrand looks particularly simple.
Let $x_k$ be one of the $|W|$ fixed points of $g_0$. Then, for every
$g \in \Omega$, the point $\omega(g) \cdot x_k$ is a fixed point
of $g$.
As before, let $\g b_k \subset \g g$ be the Borel subalgebra
determined by $x_k$: $\g b_k$ consists of all elements of
$\g g$ fixing $x_k$. This Borel subalgebra $\g b_k$ in turn
determines a positive root system
such that $\g b_k$ contains $\g t$ and all the negative root spaces.
Let $\g n_k$ be the nilpotent subalgebra of $\g g$ containing all
the positive root spaces with respect to $\g b_k$, so that
$\g g = \g b_k \oplus \g n_k$ as linear spaces.
We define a map
$\psi_k: \Omega \times \g n_k \to \Omega \times X$,
$$
(g, n) \mapsto (g, \omega(g) \exp(n) \cdot x_k).
$$
Then each $\psi_k$ is a diffeomorphism onto its image, and their images
for $k=1, \dots, |W|$ cover all of $\Omega \times X$.
Thus they form an atlas $\{\psi_1,\dots,\psi_{|W|}\}$ of $\Omega \times X$.

Also, for each $g \in \Omega$, we define a map $\psi_{g,k}: \g n_k \to X$,
$$
\psi_{g,k}(n)=\psi_k(g,n)=\omega(g)\exp(n) \cdot x_k.
$$
Then, for each $g \in \Omega$, $\{\psi_{g,1},\dots,\psi_{g,|W|}\}$
is an atlas of $X$. Notice that when $g=g_0$ or, more generally,
$g \in \Omega \cap \g t_{\BB R}$, $\omega(g) = e$ and the maps
$\psi_{g,1},\dots,\psi_{g,|W|}$ are the same as the maps
$\psi_{g_0,1},\dots,\psi_{g_0,|W|}$ defined at the
beginning of the previous section.

We decompose $\g n_k$ into root spaces with respect to $\g t$:
$\g n_k= \bigoplus_{l=1}^n \g n_{\alpha_{x_k,l}}$.
Each $\g n_{\alpha_{x_k,l}}$ is a one-dimensional complex space, so we can
choose a linear coordinate $z_l: \g n_{\alpha_{x_k,l}} \tilde \to \BB C$.
Let $\tilde g = \omega(g)^{-1} g \omega(g)$,
$\tilde g \in \g t_{\BB R}$.
Then by Lemma \ref{linearvf}, the vector field
$\operatorname{VF}_{\tilde g}$ on $X$ generated by
$\tilde g$ is expressed in the coordinate system
$\psi_{\tilde g,k} \equiv \psi_{g_0,k}$ by
$$
\alpha_{x_k,1}(\tilde g) z_1 \frac{\partial}{\partial z_1} + \dots +
\alpha_{x_k,n}(\tilde g) z_n \frac{\partial}{\partial z_n}.
$$
On the other hand, $\operatorname{VF}_g$ has the same
expression in the coordinate system $\psi_{g,k}$ as
$\operatorname{VF}_{\tilde g}$ does in $\psi_{\tilde g,k}$.
Hence the vector field $\operatorname{VF}_g$ on $X$ generated by
$g$ is expressed in the coordinate system $\psi_k$ by
\begin{multline}  \label{vfield'}
\alpha_{x_k,1}(\tilde g) z_1 \frac{\partial}{\partial z_1} + \dots +
\alpha_{x_k,n}(\tilde g) z_n \frac{\partial}{\partial z_n}   \\
= \alpha_{x_k,1}(g) z_1 \frac{\partial}{\partial z_1} + \dots +
\alpha_{x_k,n}(g) z_n \frac{\partial}{\partial z_n}.
\end{multline}

Expand $z_1, \dots , z_n$ to a standard coordinate system
$z_1, \dots, z_n,\xi_1, \dots , \xi_n$ on $T^* \g n_k$ so that
every element of
$T^* \g n_k \simeq \g n_k \times \g n_k^*$
is expressed in these coordinates as
$(z_1,\dots, z_n, \xi_1 dz_1 + \dots + \xi_n dz_n)$.
This gives us a chart
$$
\tilde \psi_k: (g, z_1,\dots,z_n,\xi_1,\dots,\xi_n)
\to \Omega \times T^*X
$$
and an atlas $\{ \tilde \psi_1,\dots,\tilde \psi_{|W|} \}$
of $\Omega \times T^*X$.

Because $G_{\BB R}$ acts on $X$ by complex automorphisms,
the differential form $\sigma$ in these coordinates is
$d\xi_1 \wedge dz_1 + \dots +d\xi_n \wedge dz_n$.

Recall that $\mu_{\lambda}=\mu+\lambda_x$, hence the integrand
$e^{\langle g, \mu_{\lambda}(\zeta) \rangle} \phi
(-\sigma + \pi^* \tau_{\lambda})^n$ becomes
$$
n! e^{\langle g, \lambda_x \rangle}
\exp \bigl(\alpha_{x_k,1}(g)z_1\xi_1 + \dots + \alpha_{x_k,n}(g)z_n\xi_n \bigr)
\phi(g) (-\sigma + \pi^* \tau_{\lambda})^n.
$$

Let $\|z\|_k= \sqrt{|z_1|^2 + \dots + |z_n|^2}$ and
$\|\xi\|_k=\sqrt{|\xi_1|^2 + \dots + |\xi_n|^2}$; these norms
are defined inside respective charts and not globally on $T^*X$.

\separate

Find $D$ sufficiently large so that the disks around fixed points
of $g_0$
$$\{\psi_1(g_0,z);\,\|z\|_1 \le D\}, \dots,
\{\psi_{|W|}(g_0,z);\,\|z\|_{|W|} \le D\}
$$
cover all of $X$.

Next find an $\epsilon>0$ small enough so that for each $k=1,\dots,|W|$
\begin{equation} \label{one}
\{\psi_k(g_0,z);\,\|z\|_k \le \epsilon\} \cap
\bigcup_{l \ne k} \{\psi_l(g_0,z);\,\|z\|_l \le 3D\} = \varnothing.
\end{equation}
We will also assume that $\epsilon \le D/2$.

Let $\delta: \BB R \to [0,1]$ be a smooth bump function which
takes on value $1$ on $[-D, D]$, vanishes outside
$(-2D,2D)$, and is  nondecreasing
on negative reals, nonincreasing on positive reals.

Let $\gamma: \BB R^+ \to (0,1]$ be another smooth function
such that $\gamma([0,1])=\{1\}$,
$\gamma(x)=\frac 1x$ for $x >2$,
$\frac 1x \le \gamma(x) \le \frac 2x$ for all $x \ge 1$,
and $\gamma$ is nonincreasing.

\separate

For each $t \in [0,1]$ and $k=1,\dots,|W|$, we will define a map
$\Theta_t^k: \Omega \times T^*X \to \Omega \times T^*X$.
First of all we define a diffeomorphism $\tilde \Theta_t^k$
on $T^* \g n_k \simeq \g n_k \times \g n_k^*$ by
\begin{multline*}
(z_1, \dots, z_n,\xi_1, \dots, \xi_n) \mapsto \\
\bigl( z_1 - \frac {\bar \alpha_{x_k,1}(g)}{|\alpha_{x_k,1}(g)|}
t \epsilon \delta(\|z\|_k)\gamma(t \|\xi\|_k) \bar \xi_1, \dots,
z_n - \frac{\bar \alpha_{x_k,n}(g)}{|\alpha_{x_k,n}(g)|}
t \epsilon \delta(\|z\|_k)\gamma(t \|\xi\|_k) \bar \xi_n, \\
\xi_1, \dots, \xi_n \bigr).
\end{multline*}

$\tilde \Theta_t^k$ shifts $(z_1,\dots,z_n)$ by a vector
$$
- t \epsilon \delta(\|z\|_k)\gamma(t \|\xi\|_k)
\bigl( \frac {\bar \alpha_{x_k,1}(g)}{|\alpha_{x_k,1}(g)|} \bar \xi_1, \dots,
\frac {\bar \alpha_{x_k,n}(g)}{|\alpha_{x_k,n}(g)|} \bar \xi_n \bigr),
$$
which has length at most $2\epsilon \le D$
(because $\gamma(x) \le \frac 2x$).
Hence the maps $\tilde \Theta^k_t$ and $(\tilde \Theta_t^k)^{-1}$
leave points outside the set $\{(z,\xi);\,\|z\|_k \le 3D\}$
completely unaffected.
Then we use the diffeomorphism between $\g n_k \times \g n_k^*$ and
$T^*(\exp(\g n_k) \cdot x_k) \subset T^*X$ induced by the map
$\psi_{g_0,k}: \g n_k \to X$, $\psi_{g_0,k}(n) = \exp(n) \cdot x_k$,
to regard $\tilde \Theta_t^k$ as a map on
$T^*(\exp(\g n_k) \cdot x_k)$.
But since $\tilde \Theta_t^k$ becomes the identity map when the
basepoint of $\zeta \in T^*X$ lies away from the compact subset
$$
\psi_{g_0,k} (\{z;\, \|z\|_k \le 3D \}) \subset
\exp(\g n_k) \cdot x_k \subset X,
$$
$\tilde \Theta_t^k$
can be extended by identity to a diffeomorphism $T^*X \to T^*X$.

Finally, we define
$\Theta_t^k: \Omega \times T^*X \to \Omega \times T^*X$
using the ``twisted'' product structure of $\Omega \times T^*X$
induced by $\omega(g)$.
Recall that the group $G$ acts on $X$ which induces an action on $T^*X$.
For $\tilde \gamma \in G$ and $\zeta \in T^*X$, we denote this action by
$\tilde \gamma \cdot \zeta$.
Then, for $(g,\zeta) \in \Omega \times T^*X$, we set
$$
\Theta_t^k (g,\zeta) = \bigl(
g, \omega(g) \cdot (\tilde \Theta_t^k (\omega(g)^{-1} \cdot \zeta)) \bigr).
$$

\separate

Inside the chart $\tilde \psi_k$ centered at the point $(g_0,x_k)$,
$\Theta_t^k$ is formally given by the same expression as before:
\begin{multline*}
(g, z_1, \dots, z_n,\xi_1, \dots, \xi_n) \mapsto \\
\bigl( g, z_1 - \frac {\bar \alpha_{x_k,1}(g)}{|\alpha_{x_k,1}(g)|}
t \epsilon \delta(\|z\|_k)\gamma(t \|\xi\|_k) \bar \xi_1, \dots,
z_n - \frac {\bar \alpha_{x_k,n}(g)}{|\alpha_{x_k,n}(g)|}
t \epsilon \delta(\|z\|_k)\gamma(t \|\xi\|_k) \bar \xi_n,  \\
\xi_1, \dots, \xi_n \bigr),
\end{multline*}
that is we shift $(z_1,\dots,z_n)$ by a vector
\begin{equation}   \label{zshift}
- t \epsilon \delta(\|z\|_k)\gamma(t \|\xi\|_k)
\bigl( \frac {\bar \alpha_{x_k,1}(g)}{|\alpha_{x_k,1}(g)|} \bar \xi_1, \dots,
\frac {\bar \alpha_{x_k,n}(g)}{|\alpha_{x_k,n}(g)|} \bar \xi_n \bigr).
\end{equation}

This choice of coefficients
$- \frac {\bar \alpha_{x_k,l}(g)}{|\alpha_{x_k,l}(g)|}$
and the equations (\ref{mu}), (\ref{vfield'}) imply that
\begin{equation}  \label{realpart}
Re((\Theta_t^k)^*\langle g, \mu(\zeta) \rangle)
\le  Re(\langle g, \mu(\zeta) \rangle),
\end{equation}
and the equality occurs if and only if
$\Theta_t^k(g,\zeta)=(g,\zeta)$.

\separate

We define $\Theta_t: \Omega \times T^*X \to \Omega \times T^*X$ by
$$
\Theta_t= \Theta_t^{|W|} \circ \dots \circ \Theta_t^1.
$$
Observe that $\Theta_0$ is the identity map.
The following four lemmas are some of the key properties of $\Theta_t$
that we will use.

\separate

\begin{lem}   \label{equal}
For each $k=1,\dots,|W|$, the maps $\Theta_t$ and $\Theta_t^k$ coincide
on the set $\{\tilde \psi_k(g,z,\xi);\,g \in \Omega, \|z\|_k \le \epsilon\}
\subset \Omega \times T^*X$.
\end{lem}

\pf Follows immediately from condition (\ref{one}).
\qed

\begin{lem} \label{realpartlemma}
If $t>0$ and $\zeta \in T^*X$ does not lie in the zero section,
$$
Re \bigl( (\Theta_t)^*\langle g, \mu(\zeta) \rangle \bigr)
<  Re(\langle g, \mu(\zeta) \rangle).
$$
\end{lem}

\pf By (\ref{realpart}), we have
$$
Re \bigl( (\Theta_t)^*\langle g, \mu(\zeta) \rangle \bigr)
\le  Re(\langle g, \mu(\zeta) \rangle),
$$
and the equality is possible only if $\Theta^k_t(g,\zeta)=(g,\zeta)$
for all $k=1,\dots,|W|$. Because of our choice of $D$, it means that
the equality is possible only if $t=0$ or $\zeta$ lies in the zero section.
\qed

\begin{lem} \label{scaling}
There exists an $R_0>0$ (depending on $t$)
such that whenever $g \in \operatorname{supp}(\phi)$,
$\zeta \in T^*X$ and $\|\zeta\| \ge R_0$ we have
$\Theta_t (g,E \zeta)=E \Theta_t (g,\zeta)$ for all real $E \ge 1$.
That is $\Theta_t$ almost commutes with scaling the fiber.

Moreover, there is an $\tilde R_0>0$, independent of $t \in (0,1]$,
such that $R_0$ can be chosen to be $\tilde R_0 /t$.
\end{lem}

\pf
Recall that $\Theta_t= \Theta_t^{|W|} \circ \dots \circ \Theta_t^1$.
We will prove by induction on $k$ that there exists an $\tilde R_0>0$
such that whenever $g \in \operatorname{supp}(\phi)$
and $\|\zeta\| \ge \tilde R_0/t$,
$$
(\Theta_t^k \circ \dots \circ \Theta_t^1) (g,E \zeta)=
E (\Theta_t^k \circ \dots \circ \Theta_t^1) (g,\zeta)
$$
for all real $E \ge 1$.

When $\|\xi\|_k \ge 2/t$, $\gamma(t \|\xi\|_k) = \frac 1{t \|\xi\|_k}$
and the shift vector in (\ref{zshift})
\begin{multline}   \label{signs}
- t \epsilon \delta(\|z\|_k)\gamma(t \|\xi\|_k)
\bigl( \frac {\bar \alpha_{x_k,1}(g)}{|\alpha_{x_k,1}(g)|} \bar \xi_1, \dots,
\frac {\bar \alpha_{x_k,n}(g)}{|\alpha_{x_k,n}(g)|} \bar \xi_n \bigr)  \\
= - \frac{\epsilon\delta(\|z\|_k)}{\|\xi\|_k}
\bigl( \frac {\bar \alpha_{x_k,1}(g)}{|\alpha_{x_k,1}(g)|} \bar \xi_1, \dots,
\frac {\bar \alpha_n(g)}{|\alpha_n(g)|} \bar \xi_n \bigr)
\end{multline}
stays unchanged if we replace $(\xi_1,\dots,\xi_n)$ with
$(E\xi_1,\dots,E\xi_n)$, for any real $E \ge 1$.
Hence in this situation
$\Theta_t^k(g,E\zeta) = E \Theta_t^k(g,\zeta)$.

By induction hypothesis there is an $\tilde R_0>0$ such that
whenever $g \in \operatorname{supp} (\phi)$ and
$\|\zeta\| \ge \tilde R_0/t$,
$$
(\Theta_t^{k-1} \circ \dots \circ \Theta_t^1) (g,E \zeta)=
E (\Theta_t^{k-1} \circ \dots \circ \Theta_t^1) (g,\zeta).
$$

Making $\tilde R_0$ larger if necessary, we can assume in addition
that whenever
$\|z((\Theta_t^{k-1} \circ \dots \circ \Theta_t^1) (g,\zeta))\|_k \le 2D$
and $\|\zeta\| \ge \tilde R_0$ we have
$\|\xi((\Theta_t^{k-1} \circ \dots \circ \Theta_t^1) (g,\zeta))\|_k \ge 2$,
for all $t \in [0,1]$.

Then, for $t \in (0,1]$,
$\|z((\Theta_t^{k-1} \circ \dots \circ \Theta_t^1) (g,\zeta))\|_k \le 2D$
and $\|\zeta\| \ge \tilde R_0/t$ imply
$\|\xi((\Theta_t^{k-1} \circ \dots \circ \Theta_t^1) (g,\zeta))\|_k \ge 2/t$,
and so
$$
(\Theta_t^k \circ \dots \circ \Theta_t^1) (g,E \zeta)=
E (\Theta_t^k \circ \dots \circ \Theta_t^1) (g,\zeta).
$$

On the other hand, if 
$\|z((\Theta_t^{k-1} \circ \dots \circ \Theta_t^1) (g,\zeta))\|_k > 2D$,
then by induction hypothesis $\|\zeta\| \ge \tilde R_0/t$ implies
\begin{multline*}
(\Theta_t^k \circ \dots \circ \Theta_t^1) (g,E \zeta)
= (\Theta_t^{k-1} \circ \dots \circ \Theta_t^1) (g,E \zeta)  \\
= E (\Theta_t^{k-1} \circ \dots \circ \Theta_t^1) (g,\zeta)
= E (\Theta_t^k \circ \dots \circ \Theta_t^1) (g,\zeta).
\qed
\end{multline*}

Now let us return to calculation of our integral (\ref{theintegral}).
Instead of integrating our original form
$e^{\langle g, \mu_{\lambda}(\zeta) \rangle} \phi
(-\sigma + \pi^* \tau_{\lambda})^n$
we fix $t \in (0,1)$ and integrate
$\Theta_t^*(e^{\langle g, \mu_{\lambda}(\zeta) \rangle} \phi
(-\sigma + \pi^* \tau_{\lambda})^n)$.
Doing so is justified by the following lemma which will be
proved in Section \ref{proofslanting}.

\begin{lem} \label{slanting}
For any $t \in [0,1]$, we have:
\begin{multline*}
\lim_{R \to \infty}
\int_{\g g_{\BB R}' \times (Ch({\cal F}) \cap \{\|\zeta\| \le R\})}
\bigl( e^{\langle g, \mu_{\lambda}(\zeta) \rangle} \phi
(-\sigma + \pi^* \tau_{\lambda})^n   \\
- \Theta_t^*(e^{\langle g, \mu_{\lambda}(\zeta) \rangle} \phi
(-\sigma + \pi^* \tau_{\lambda})^n) \bigr) =0.
\end{multline*}
\end{lem}

Recall the Borel-Moore chain $C(g_0)$ described in Proposition \ref{C(g)prop}.
Because $\Omega$ was chosen so that both
$\Omega$ and $\Omega \cap \g t_{\BB R}$ are connected,
by Remark \ref{C(g)remark}, for each $g \in \Omega \cap \g t_{\BB R}$,
we can choose $C(g)$ equal $C(g_0)$.
Moreover, for each $g \in \Omega$, we can choose $C(g)$ equal
$\omega(g)_*C(g_0)$. These chains $C(g)$, $g \in \Omega$, piece together
into a Borel-Moore chain in $\Omega \times T^*X$
of dimension $\dim_{\BB R} \g g_{\BB R}+2n+1$
which appears in each chart $\tilde \psi_k$ as $\Omega \times C(g_0)$,
$$
\partial C = \Omega \times Ch({\cal F}) -
\Omega \times (m_1 T^*_{\omega(g) \cdot x_1}X + \dots
+ m_{|W|} T^*_{\omega(g) \cdot x_{|W|}}X)
$$
and the support of $C$ lies inside
$\{(g,\zeta) \in \Omega \times T^*X;\,
Re( \langle g,\mu(\zeta) \rangle) \le 0\}$.

Take an $R \ge 1$ and restrict all cycles to the set
$\{ (g,\zeta) \in \Omega \times T^*X;\, \|\zeta\| \le R\}$.
Let $C_{\le R}$ denote the restriction of the cycle $C$, then
it has boundary
\begin{multline*}
\partial C_{\le R} =
\Omega \times (Ch({\cal F}) \cap \{\|\zeta\| \le R\}) - C'(R) \\
- \Omega \times
\bigl( m_1 (T^*_{\omega(g) \cdot x_1}X \cap \{\|\zeta\| \le R\}) + \dots
+ m_{|W|} (T^*_{\omega(g) \cdot x_{|W|}}X \cap \{\|\zeta\| \le R\}) \bigr),
\end{multline*}
where $C'(R)$ is a $(\dim_{\BB R} \g g_{\BB R} +2n)$-chain supported in the set
$$
\{ (g, \zeta) \in \Omega \times T^*X;\, \|\zeta\|=R,\,
Re( \langle g, \mu(\zeta) \rangle ) \le 0 \}.
$$
Because the chain $C$ is conic, the piece of boundary $C'(R)$
depends on $R$ by appropriate scaling of $C'(1)$ in the fiber direction.

\begin{lem}
For a fixed $t \in (0,1]$,
$$
\lim_{R \to \infty} \int_{C'(R)} 
\Theta_t^* \bigl( e^{\langle g, \mu_{\lambda}(\zeta) \rangle} \phi
(-\sigma + \pi^* \tau_{\lambda})^n \bigr) =0.
$$
\end{lem}

\pf
Let $R_0$ be as in Lemma \ref{scaling}.
Then, for $R \ge R_0$, the chain
$(\Theta_t)_* C'(R)$ depends on $R$ by scaling
$(\Theta_t)_* C'(R_0)$ in the fiber direction.

Recall that the coefficients
$- \frac {\bar \alpha_{x_k,l}(g)}{|\alpha_{x_k,l}(g)|}$
in (\ref{zshift}) are chosen so that for each
$l=1,\dots,n$ the term
$\alpha_{x_k,l}(g) \xi_l
\bigl( - \frac {\bar \alpha_{x_k,l}(g)}{|\alpha_{x_k,l}(g)|} \bar \xi_l \bigr)
= - |\alpha_{x_k,l}(g)||\xi_l|^2$
is negative.
Hence, for every $(g, \zeta)$ lying in the support of
$C'(R)$, the real part of
$\langle g,\mu(\Theta_t(g, \zeta)) \rangle$
is strictly negative.
By compactness of
$|(\Theta_t)_*C'(R_0)| \cap (\supp(\phi) \times T^*X)$,
there exists an $\epsilon'>0$ such that, whenever $(g, \zeta)$
lies in the support of $(\Theta_t)_*C'(R_0)$ and $g$ lies
in the support of $\phi$, we have
$Re( \langle g,\mu(\zeta) \rangle ) \le -\epsilon'$.
Then, for all $R \ge R_0$ and all
$(g, \zeta) \in
|(\Theta_t)_*C'(R)| \cap (\supp(\phi) \times T^*X)$, we have
$Re( \langle g,\mu(\zeta) \rangle ) \le -\epsilon' \frac R{R_0}$.

Integrating the form
$\Theta_t^* \bigl( e^{\langle g, \mu_{\lambda}(\zeta) \rangle} \phi
(-\sigma + \pi^* \tau_{\lambda})^n \bigr)$
over the chain $C'(R)$ is equivalent to integrating
$e^{\langle g, \mu_{\lambda}(\zeta) \rangle} \phi
(-\sigma + \pi^* \tau_{\lambda})^n =
e^{\langle g, \mu (\zeta) \rangle}
e^{\langle g, \lambda_x (\zeta) \rangle}
\phi (-\sigma + \pi^* \tau_{\lambda})^n$
over $(\Theta_t)_* C'(R)$.
Since the integrand decays exponentially over the
support of $(\Theta_t)_* C'(R)$, the integral
tends to zero as $R \to \infty$.
\qed

Thus
\begin{multline*}
\int_{Ch({\cal F})}
\mu_{\lambda}^* \hat \phi (-\sigma+\pi^* \tau_{\lambda})^n  \\
=\lim_{R \to \infty}
\int_{\Omega \times (Ch({\cal F}) \cap \{\|\zeta\| \le R\})}
e^{\langle g, \mu_{\lambda}(\zeta) \rangle} \phi
(-\sigma+\pi^* \tau_{\lambda})^n  \\
=\lim_{R \to \infty}
\int_{\Omega \times (Ch({\cal F}) \cap \{\|\zeta\| \le R\})}
\Theta_t^*(e^{\langle g, \mu_{\lambda}(\zeta) \rangle} \phi
(-\sigma + \pi^* \tau_{\lambda})^n)  \\
=\lim_{R \to \infty}
\int_{C'(R) + \Omega \times \bigl(
\sum_{k=1}^{|W|} m_k (T^*_{\omega(g) \cdot x_k}X \cap \{\|\zeta\| \le R\})
\bigr)}
\Theta_t^*(e^{\langle g, \mu_{\lambda}(\zeta) \rangle} \phi
(-\sigma + \pi^* \tau_{\lambda})^n)  \\
=\lim_{R \to \infty}
\int_{\Omega \times \bigl(
\sum_{k=1}^{|W|} m_k (T^*_{\omega(g) \cdot x_k}X \cap \{\|\zeta\| \le R\})
\bigr)}
\Theta_t^*(e^{\langle g, \mu_{\lambda}(\zeta) \rangle} \phi
(-\sigma + \pi^* \tau_{\lambda})^n),
\end{multline*}
i.e. the integral over $C'(R)$ can be ignored and
we are left with integrals over
$m_k \bigl( \Omega \times
(T^*_{\omega(g) \cdot x_k}X \cap \{\|\zeta\| \le R\} ) \bigr)$,
for $k=1,\dots,|W|$.
Because the integral converges absolutely, we can let $R \to \infty$ and
drop the restriction $\|\zeta\|\le R$:
\begin{equation}  \label{integral2}
\int_{Ch({\cal F})}
\mu_{\lambda}^* \hat \phi (-\sigma+\pi^* \tau_{\lambda})^n =
\int_{\Omega \times \bigl(
\sum_{k=1}^{|W|} m_k T^*_{\omega(g) \cdot x_k}X \bigr)}
\Theta_t^*(e^{\langle g, \mu_{\lambda}(\zeta) \rangle} \phi
(-\sigma + \pi^* \tau_{\lambda})^n).
\end{equation}

Lemma \ref{equal} tells us that the maps
$\Theta_t$ and $\Theta_t^k$ coincide over $T^*_{\omega(g) \cdot x_k}X$:
$$
\Theta_t|_{T^*_{\omega(g) \cdot x_k}X} \equiv
\Theta_t^k|_{T^*_{\omega(g) \cdot x_k}X}.
$$
We also have $\delta(\|z\|_k)=1$ and
the exponential part
$\Theta_t^* \bigl( \langle g, \mu_{\lambda}(\zeta) \rangle \bigl)$
of our integrand
$\Theta_t^* \bigl( e^{\langle g, \mu_{\lambda}(\zeta) \rangle} \phi
(-\sigma + \pi^* \tau_{\lambda})^n \bigr)$
becomes
\begin{equation}  \label{exp}
- t \epsilon \gamma(t\|\xi\|_k) \bigl( |\alpha_1(g)| \xi_1\bar \xi_1 +\dots+
|\alpha_n(g)| \xi_n \bar \xi_n \bigr) +
\Theta_t^* \langle g,\lambda_x \rangle.
\end{equation}

We know that
$\int_{\g g_{\BB R}' \times Ch({\cal F})}
\Theta_t^* \bigl( e^{\langle g, \mu_{\lambda}(\zeta) \rangle} \phi
(-\sigma + \pi^* \tau_{\lambda})^n \bigr)$
does not depend on $t$. So in order to calculate its value
we are allowed to regard it as a constant function of $t$
and take its limit as $t \to 0^+$.

We can break up our chain $m_k(T^*_{\omega(g) \cdot x_k}X)$ into
two portions: one portion where $\|\xi(g,\zeta)\|_k \ge 1/t$ and the
other where $\|\xi(g,\zeta)\|_k < 1/t$.

\begin{lem}
$$
\lim_{t \to 0^+}
\int_{m_k \bigl( \Omega \times
(T^*_{\omega(g) \cdot x_k}X \cap \{\|\xi(g,\zeta)\|_k \ge 1/t\}) \bigr)}
\Theta_t^* \bigl( e^{\langle g, \mu_{\lambda}(\zeta) \rangle} \phi
(-\sigma + \pi^* \tau_{\lambda})^n \bigr) =0.
$$
\end{lem}

\pf
When $\|\xi\|_k \ge 1/t$,
$\gamma(t \|\xi\|_k) \ge \frac 1{t \|\xi\|_k}$
and the exponential part (\ref{exp}) is at most
$$
- \frac {\epsilon}{\|\xi\|_k} \bigl( |\alpha_{x_k,1}(g)| \xi_1\bar \xi_1
+\dots+ |\alpha_{x_k,n}(g)| \xi_n \bar \xi_n \bigr) +
\Theta_t^* \langle g,\lambda_x \rangle.
$$
But $\xi_1 \bar \xi_1 + \dots + \xi_n \bar \xi_n = \|\xi\|_k^2$, so
at least one of the $\xi_l \bar \xi_l \ge \|\xi\|_k^2/n$.
Thus we get a new estimate of (\ref{exp}) from above:
$$
- \frac {\epsilon}n |\alpha_{x_k,l}(g)| \|\xi\|_k
+ \Theta_t^* \langle g,\lambda_x \rangle
\le - \frac {\epsilon}{nt} |\alpha_{x_k,l}(g)|
+ \Theta_t^* \langle g,\lambda_x \rangle.
$$
Because the term
$\Theta_t^* \langle g,\lambda_x \rangle$
is bounded, the last expression tends to $-\infty$ as $t \to 0^+$,
i.e. the integrand decays exponentially and the lemma follows.
\qed

Thus, in the formula (\ref{integral2}) the integral over the portion
$$
m_k \bigl( \Omega \times (T^*_{\omega(g) \cdot x_k}X \cap
\{\|\xi(g,\zeta)\|_k \ge 1/t\}) \bigr)
$$
can be ignored too:
\begin{multline*}
\int_{Ch({\cal F})}
\mu_{\lambda}^* \hat \phi (-\sigma+\pi^* \tau_{\lambda})^n   \\
= \lim_{t \to 0^+} \int_{\sum_{k=1}^{|W|}
m_k \bigl( \Omega \times (T^*_{\omega(g) \cdot x_k}X \cap
\{\|\xi(g,\zeta)\|_k < 1/t\}) \bigr)}
\Theta_t^*(e^{\langle g, \mu_{\lambda}(\zeta) \rangle} \phi
(-\sigma + \pi^* \tau_{\lambda})^n).
\end{multline*}

Finally, over the portion
$m_k \bigl( \Omega \times (T^*_{\omega(g) \cdot x_k}X \cap
\{\|\xi(g,\zeta)\|_k < 1/t\}) \bigr)$,
the function $\gamma(t \|\xi\|_k)$ is identically one,
so the exponential part (\ref{exp}) reduces to
$$
- t \epsilon \bigl( |\alpha_{x_k,1}(g)| \xi_1\bar \xi_1 +\dots+
|\alpha_{x_k,n}(g)| \xi_n \bar \xi_n \bigr) +
\Theta_t^* \langle g,\lambda_x \rangle.
$$
We also have $\Theta_t^*(\phi)=\phi$,  $\Theta_t^*(d\xi_l)=d\xi_l$,
$$
\Theta_t^*(dz_l) = -d \bigl( t \epsilon \gamma(t\|\xi\|_k)
\frac {\bar \alpha_{x_k,l}(g)} {|\alpha_{x_k,l}(g)|} \bar \xi_l \bigr)
=-t\epsilon \frac {\bar \alpha_{x_k,l}(g)} {|\alpha_{x_k,l}(g)|} d\bar \xi_l,
$$
and, because $\tau_{\lambda}$ is a 2-from on $X$,
\begin{multline*}
\Theta_t^* (\pi^* \tau_{\lambda})   \\
= t^2 \sum_{a,b=1}^n \bigl( A_{a,b}(g,t\xi_1,\dots,t\xi_n) d\xi_a \wedge d\xi_b
+ B_{a,b}(g,t\xi_1,\dots,t\xi_n) d \bar\xi_a \wedge d\xi_b   \\
+ C_{a,b}(g,t\xi_1,\dots,t\xi_n) d\xi_a \wedge d\bar\xi_b +
D_{a,b}(g,t\xi_1,\dots,t\xi_n) d\bar\xi_a \wedge d\bar\xi_b \bigr),
\end{multline*}
where each $A_{a,b},B_{a,b},C_{a,b},D_{a,b}$ is a bounded function of
$(g,t\xi_1,\dots,t\xi_n)$.
Similarly,
$$
\Theta_t^* \langle g,\lambda_x \rangle =
\langle g,\lambda_{\omega(g) \cdot x_k} \rangle +
t \sum_{c=1}^n \bigl( \xi_c E_c(g,t\xi_1,\dots,t\xi_n) +
\bar\xi_c F_c(g,t\xi_1,\dots,t\xi_n) \bigr)
$$
for some bounded functions $E_c$,
$F_c$ of $(g,t\xi_1,\dots,t\xi_n)$, $c=1,\dots,n$.

Thus we end up integrating
\begin{multline*}
n! e^{\langle g,\lambda_{\omega(g) \cdot x_k} \rangle}
e^{-t \epsilon \bigl( |\alpha_{x_k,1}(g)| \xi_1\bar \xi_1 +\dots+
|\alpha_{x_k,n}(g)| \xi_n \bar \xi_n \bigr)}   \\
\cdot e^{t \sum_{c=1}^n \bigl( \xi_c E_c(g,t\xi_1,\dots,t\xi_n) +
\bar\xi_c F_c(g,t\xi_1,\dots,t\xi_n) \bigr)}   \\
\cdot \phi(g) \wedge
\Bigl( (t \epsilon)^n \frac {\bar \alpha_{x_k,1}(g)} {|\alpha_{x_k,1}(g)|}
\dots \frac {\bar \alpha_{x_k,n}(g)} {|\alpha_{x_k,n}(g)|}
d\xi_1 \wedge d\bar\xi_1 \wedge \dots \wedge d\xi_n \wedge d\bar\xi_n   \\
+ \text{terms containing $\Theta_t^* (\pi^* \tau_{\lambda})$} \Bigr)
\end{multline*}
over $m_k \bigl( \Omega \times (T^*_{\omega(g) \cdot x_k}X \cap
\{\|\xi(g,\zeta)\|_k < 1/t\}) \bigr)$.
(Recall that the orientation of this chain is determined by the
product orientation on $\Omega \times T^*_{\omega(g) \cdot x_k}X$,
and the orientation of $T^*_{\omega(g) \cdot x_k}X$ is given by
(\ref{orientation}).)
Changing variables $y_l = \sqrt{\epsilon t} \xi_l$ for $l=1,\dots,n$
we obtain
\begin{multline*}
m_k n!
\int_{\Omega} e^{\langle g,\lambda_{\omega(g) \cdot x_k} \rangle} \phi(g)
\int_{\{|y_1|^2+\dots+|y_n|^2 < \frac {\epsilon}t \}}
e^{-|\alpha_{x_k,1}(g)| |y_1|^2 -\dots- |\alpha_{x_k,n}(g)| |y_n|^2}    \\
\cdot e^{\sqrt{t} \sum_{c=1}^n
\bigl( y_c E_c(g,\sqrt{t}y_1,\dots,\sqrt{t}y_n) +
\bar y_c F_c(g,\sqrt{t}y_1,\dots,\sqrt{t}y_n) \bigr)}    \\
\cdot \Bigl( \frac {\bar \alpha_{x_k,1}(g)} {|\alpha_{x_k,1}(g)|} \dots
\frac {\bar \alpha_{x_k,n}(g)} {|\alpha_{x_k,n}(g)|}
dy_1 \wedge d\bar y_1 \wedge \dots \wedge dy_n \wedge d\bar y_n
+ t \cdot (\text{bounded terms}) \Bigr).
\end{multline*}
Because the term
$$
\sqrt{t} \sum_{c=1}^n
\bigl( y_c E_c(g,\sqrt{t}y_1,\dots,\sqrt{t}y_n) +
\bar y_c F_c(g,\sqrt{t}y_1,\dots,\sqrt{t}y_n) \bigr)
$$
can be bounded on $\{|y_1|^2+\dots+|y_n|^2 < \frac {\epsilon}t \}$
independently of $t$,
by the Lebesgue dominant convergence theorem this integral tends to
\begin{multline*}
m_k n!
\int_{\Omega} e^{\langle g,\lambda_{\omega(g) \cdot x_k} \rangle} \phi(g)
\int_{\{(y_1, \dots, y_n) \in \BB C^n\}}
e^{-|\alpha_{x_k,1}(g)| |y_1|^2 -\dots- |\alpha_{x_k,n}(g)| |y_n|^2}    \\
\cdot \frac {\bar \alpha_{x_k,1}(g)} {|\alpha_{x_k,1}(g)|} \dots
\frac {\bar \alpha_{x_k,n}(g)} {|\alpha_{x_k,n}(g)|}
dy_1 \wedge d\bar y_1 \wedge \dots \wedge dy_n \wedge d\bar y_n  \\
= m_k n! (2\pi i)^n \int_{\Omega}
\frac {e^{\langle g,\lambda_{\omega(g) \cdot x_k} \rangle} \phi(g)}
{\alpha_{x_k,1}(g) \dots \alpha_{x_k,n}(g)}
\end{multline*}
as $t \to 0^+$. Therefore, this is the value of our original integral
(\ref{integral_formula}).

The last expression may appear to be missing a factor of $(-1)^n$, but it is
correct because we are using the convention (11.1.2) in \cite{KaScha}
to identify the real cotangent bundle $T^*(X^{\BB R})$ with the
holomorphic cotangent bundle $T^*X$ of the complex manifold $X$.
Under this identification, the standard symplectic
structure on $T^*(X^{\BB R})$ equals $2 Re (\sigma)$.
Moreover,
\begin{multline*}
\text{the orientation of $T^*_{\omega(g) \cdot x_k}X$ given by
(\ref{orientation})}  \\
= (-1)^n \text{ the complex orientation of $T^*_{\omega(g) \cdot x_k}X$}.
\end{multline*}

\separate

So far we assumed that the support of $\phi$ lies in some
open subset $\Omega \subset \g g_{\BB R}'$ of a special kind.
When $g$ lied inside $\Omega$ the set of points in $X$ fixed by $g$
could be expressed as
$\{ \omega(g) \cdot x_1,\dots,\omega(g) \cdot x_{|W|} \}$.
Of course, the function $\omega$ was defined on $\Omega$ only.
Now we are interested in the global situation, so,
for each $g \in \g g_{\BB R}'$, we denote by
$\{ x_1(g),\dots,x_{|W|}(g)\}$ the set of points in $X$ fixed by $g$
enumerated in a completely arbitrary way.
Similarly, when $g$ lied inside $\Omega$, each coefficient $m_k$ was constant.
Now we need to recognize that the coefficients $m_k$ are not
constant on $\g g_{\BB R}'$ globally -- they are only locally constant.
To emphasize this dependence on $g$ and the choice of fixed point
$x_k(g)$ we will write $m_{x_k(g)}$ for the multiplicity of the cycle
$T^*_{x_k(g)}X$.

We define a function $\tilde F$ on $\g g_{\BB R}'$ by
$$
\tilde F(g)= \sum_{k=1}^{|W|} 
\frac {m_{x_k(g)} e^{\langle g,\lambda_{x_k(g)}\rangle}}
{\alpha_{x_k(g),1}(g) \dots \alpha_{x_k(g),n}(g)}.
$$
A simple partition of unity argument proves the formula
\begin{equation}  \label{formula}
\frac 1{(2\pi i)^nn!} \int_{Ch({\cal F})}
\mu_{\lambda}^* \hat \phi (-\sigma+\pi^* \tau_{\lambda})^n
= \int_{\g g_{\BB R}} \tilde F \phi,
\end{equation}
when the form $\phi$ is compactly supported in $\g g_{\BB R}'$,
an open subset of the set of regular semisimple elements in
$\g g_{\BB R}$ whose complement has measure zero.
Now suppose that the support of $\phi$ does not lie inside
$\g g_{\BB R}'$.
We know from the fixed point character formula in \cite{SchV2}
that the function $\tilde F$ defined on $\g g_{\BB R}'$
is a character of a certain admissible representation of finite length,
hence by Theorem 3.3 of \cite{A} a locally $L^1$ function on $\g g_{\BB R}$.
Let $\{\phi_l\}_{l=1}^{\infty}$ be a partition of unity
on $\g g_{\BB R}'$ subordinate to the covering by those open sets
$\Omega$'s.
Then $\phi$ can be realized on $\g g_{\BB R}'$
as a pointwise convergent series:
$$
\phi= \sum_{l=1}^{\infty} \phi_l \phi.
$$
Because $\tilde F \in L^1_{loc} (\g g_{\BB R})$, the series
$\sum_{l=1}^{\infty}\int_{\g g_{\BB R}} \tilde F \phi_l \phi$
converges absolutely. Hence
\begin{multline*}
\frac 1{(2\pi i)^nn!} \int_{Ch({\cal F})}
\mu_{\lambda}^* \hat \phi (-\sigma+\pi^* \tau_{\lambda})^n   \\
= \sum_{l=1}^{\infty} \frac 1{(2\pi i)^nn!} \int_{Ch({\cal F})}
\mu_{\lambda}^* \widehat{\phi_l  \phi} (-\sigma+\pi^* \tau_{\lambda})^n   \\
= \sum_{l=1}^{\infty}\int_{\g g_{\BB R}} \tilde F \phi_l \phi =
\int_{\g g_{\BB R}} \tilde F \phi,
\end{multline*}
which proves the integral character formula (\ref{formula}) in general.

\separate

It remains to justify integrating the form
$\Theta_t^* \bigl( e^{\langle g, \mu_{\lambda}(\zeta) \rangle} \phi
(-\sigma + \pi^* \tau_{\lambda})^n \bigr)$
instead of
$e^{\langle g, \mu_{\lambda}(\zeta) \rangle} \phi
(-\sigma + \pi^* \tau_{\lambda})^n$.

\separate

\begin{section}
{Proof of Lemma \ref{slanting}}   \label{proofslanting}
\end{section}

Because the form
$e^{\langle g, \mu_{\lambda}(\zeta) \rangle} \phi
(-\sigma + \pi^* \tau_{\lambda})^n$
is closed, the integral
\begin{multline*}
\int_{\g g_{\BB R}' \times (Ch({\cal F}) \cap \{\|\zeta\| \le R\})}
\bigl( e^{\langle g, \mu_{\lambda}(\zeta) \rangle} \phi
(-\sigma + \pi^* \tau_{\lambda})^n
- \Theta_t^*(e^{\langle g, \mu_{\lambda}(\zeta) \rangle} \phi
(-\sigma + \pi^* \tau_{\lambda})^n) \bigr) \\
=\int_{\g g_{\BB R}' \times (Ch({\cal F}) \cap \{\|\zeta\| \le R\}) -
(\Theta_t)_* \bigl( \g g_{\BB R}' \times (Ch({\cal F}) \cap
\{\|\zeta\| \le R\}) \bigr)}
e^{\langle g, \mu_{\lambda}(\zeta) \rangle} \phi
(-\sigma + \pi^* \tau_{\lambda})^n
\end{multline*}
is equal to the integral of
$e^{\langle g, \mu_{\lambda}(\zeta) \rangle} \phi
(-\sigma + \pi^* \tau_{\lambda})^n$
over the chain traced by
$\Theta_{t'} \bigl( \Omega \times
\partial (Ch({\cal F}) \cap \{\|\zeta\| \le R\}) \bigr)$
as $t'$ varies from $0$ to $t$.
We will show that this integral tends to zero as $R \to \infty$.

Since $Ch({\cal F})$ is a cycle in $T^*X$, the chain
$\Omega \times \partial (Ch({\cal F}) \cap \{\|\zeta\| \le R\})$
is supported inside the set
$\Omega \times \{\zeta \in T^*X;\, \|\zeta\|=R\}$.
Because $R \to \infty$, we can assume that $R>0$.
Then the chain traced by
$\Theta_{t'} \bigl( \Omega \times
\partial (Ch({\cal F}) \cap \{\|\zeta\| \le R\}) \bigr)$
as $t'$ varies from $0$ to $t$ lies away from the zero section
$\Omega \times T^*_XX$ in $\Omega \times T^*X$.

If we regard $\Theta$ as a map
$\Omega \times T^*X \times [0,1] \to \Omega \times T^*X$,
we get an integral of
$\Theta^* \bigl( e^{\langle g, \mu_{\lambda}(\zeta) \rangle} \phi
(-\sigma + \pi^* \tau_{\lambda})^n \bigr)$
over the chain
$\Omega \times \partial (Ch({\cal F}) \cap \{\|\zeta\| \le R\}) \times [0,t]$.

The idea is to integrate out the $\Omega$ variable and check that
the result decays faster than any negative power of $R$.

Clearly, $\Theta^*(\phi)=\phi$ and Lemma \ref{realpartlemma} says that
$$
\Theta^* \langle g, \mu(\zeta) \rangle =
\langle g, \mu(\zeta) \rangle - \kappa(g,\zeta, t')
$$
for some smooth function $\kappa(g, \zeta, t')$
which has positive real part.
The integral in question can be rewritten as
$$
\int_{\Omega \times \partial (Ch({\cal F}) \cap \{\|\zeta\| \le R\})
\times [0,t]}
e^{\langle g, \mu(\zeta) \rangle}
e^{-\kappa(g, \zeta, t')} \phi(g) \Theta^*
\bigl( e^{\langle g,\lambda_x \rangle}(-\sigma+ \pi^*\tau_{\lambda})^n \bigr).
$$

We pick a system of local coordinates $(x_1,\dots,x_n)$ of $X$
and construct respective local coordinates
$(x_1,\dots,x_n,\eta_1,\dots,\eta_n)$ of $T^*X$.
Suppose that we know that all the partial
derivatives of all orders of $e^{-\kappa(g, \zeta, t')}$
and $\Theta^* (-\sigma+ \pi^*\tau_{\lambda})$
with respect to the $g$ variable can be bounded independently
of $\zeta$ and $t'$ on the set
$\supp(\phi) \times \{\zeta \in T^*X;\, \|\zeta\|>0\} \times [0,t]$.
Let $y_1, \dots, y_m$ be a system of linear coordinates on
$\g g_{\BB R}$, write
$\mu(\zeta)=\beta_1(\zeta) dy_1 + \dots +\beta_m(\zeta) dy_m$,
then
\begin{multline*}
\int_{\g g_{\BB R}}
e^{\langle g, \mu(\zeta) \rangle}
e^{-\kappa(g,\zeta, t')} \phi(g) \Theta^*
\bigl(e^{\langle g,\lambda_x \rangle}(-\sigma+ \pi^*\tau_{\lambda})^n \bigr) \\
= - \frac 1{\beta_l(\zeta)}
\int_{\g g_{\BB R}}
e^{\langle g, \mu(\zeta) \rangle}
\frac {\partial}{\partial y_l} \Bigl( e^{-\kappa(g,\zeta, t')}
\phi(g) \Theta^* \bigl(
e^{\langle g,\lambda_x \rangle}(-\sigma + \pi^*\tau_{\lambda})^n \bigr) \Bigr),
\end{multline*}
and the last integral can be bounded by a constant multiple of $R^n$.
We can keep performing integration by parts to get the desired
estimate just like for the ordinary Fourier transform.
Thus (after integrating out the $\Omega$-variable) we see that the
integrand indeed decays rapidly in the fiber variable of $T^*X$.
Hence our integral tends to zero as $R \to \infty$.

\separate

To show boundedness of the partial derivatives with respect
to $g$ we replace local coordinates
$(x_1,\dots,x_n,\eta_1,\dots,\eta_n)$ on $T^*X$ without the zero section
with ``spherical coordinates'' $(x_1,\dots,x_n,v,\nu)$, where
$$
\nu=\|(x_1,\dots,x_n,\eta_1,\dots,\eta_n)\|,
\qquad \|(x_1,\dots,x_n,v)\|=1
$$
and
$$
\nu \cdot (x_1,\dots,x_n,v) = (x_1,\dots,x_n,\eta_1,\dots,\eta_n).
$$
This change makes all variables bounded, except for only one -- $\nu$.
Recall that $y_1, \dots, y_m$ is a system of linear coordinates on
$\g g_{\BB R}$. Thus we can introduce a system of coordinates on
$\Omega \times \{\zeta \in T^*X;\, \|\zeta\|>0\} \times [0,t]$:
$$
(y_1, \dots, y_m,x_1,\dots,x_n,v,\nu,t'\nu).
$$
Notice that there is no ``twist'' by the action of $G_{\BB R}$ in
this coordinate system. It is just the direct product of the
coordinate systems for $\Omega$, $\{\zeta \in T^*X;\, \|\zeta\|>0\}$
and $[0,t]$, except that the coordinate $t'$ of $[0,t]$ got
replaced with $t'\nu$.

Because $\supp(\phi)$ is compact, the coordinates
$y_1, \dots, y_m$ can be treated as bounded.
Thus we are left with two unbounded coordinates -- $\nu$ and $t'\nu$.

A close examination of the definition of $\Theta_{t'}^k$
shows that, as long as
$$
\Theta_{t'}^k(y_1, \dots, y_m,x_1,\dots,x_n,v,\nu,t'\nu)
$$
lies in the same coordinate system,
$\Theta_{t'}^k$ can be written in these coordinates as
\begin{multline*}
\Theta_{t'}^k(y_1, \dots, y_m,x_1,\dots,x_n,v,\nu,t'\nu)  \\
= (y_1, \dots, y_m,
x_1 + \tilde x_1,\dots,x_n + \tilde x_n,v + \tilde v,\nu + \tilde \nu,
t'\nu + \widetilde{t'\nu}),
\end{multline*}
where $\tilde x_1,\dots, \tilde x_n$, $\tilde v$, $\tilde \nu$ and
$\widetilde{t'\nu}$ are some functions which
depend on $y_1, \dots, y_m$, $x_1,\dots,x_n$, $v$, $t'\nu$
and not on $\nu$. Hence same is true of $\Theta$.
On the other hand, when
$\Theta_{t'}^k(y_1, \dots, y_m,x_1,\dots,x_n,v,\nu,t'\nu)$
lies outside of this coordinate system, we can arrive to the same
conclusion about $\Theta_{t'}^k$ and $\Theta$ using another coordinate
chart of the same kind which contains the image. In other words,
$\Theta$ is different form the identity map by something that
does not depend on $\nu$.
Next we observe that Lemma \ref{scaling}
implies that $\tilde x_1,\dots, \tilde x_n$,
$\tilde v$, $\tilde \nu$ and $\widetilde{t'\nu}$
become constant with respect
to $t'\nu$ when $t'\nu$ is larger than some $\tilde R_0$.
This shows that the partial derivatives of $\Theta$ and hence
the partial derivatives of $\Theta^* (-\sigma+ \pi^*\tau_{\lambda})$
are bounded.

Note that this argument also shows that the function
$$
\frac 1{\|\zeta\|} \kappa = \frac 1{\|\zeta\|} \bigl(
\langle g, \mu(\zeta) \rangle - \Theta_{t'}^* \langle g, \mu(\zeta) \rangle
\bigr)
$$
depends on $y_1, \dots, y_m$, $x_1,\dots,x_n$, $v$, $t'\nu$
and not on $\nu$.

\separate

It remains to prove boundedness of the partial derivatives of
$e^{-\kappa(g,\zeta, t')}$ with respect to the $g$ variable.
Recall that the proof of Lemma \ref{scaling}
was based on the property of $\gamma$ that
$\gamma(t\|\xi\|_k) = \frac 1{t\|\xi\|_k}$ when
$t\|\xi\|_k > 2$.
On the other hand, the next lemma is based on the property of $\gamma$
that $\gamma(t\|\xi\|_k) = 1$ when $t\|\xi\|_k < 1$.

\separate

\begin{lem}
There exist a smooth bounded function
$\tilde \kappa (g, v, t')$ defined on
$$
\Omega \times \{\zeta \in T^*X;\, \|\zeta\|=1 \} \times [0,1]
$$
and a real number $\tilde r_0>0$
such that, whenever $t' \|\zeta\| \le \tilde r_0$,
$$
\kappa(g,\zeta, t') = t' \|\zeta\|^2 \cdot
\tilde \kappa \bigl( g, \frac{\zeta}{\|\zeta\|}, t' \bigr).
$$
Moreover, $Re(\tilde \kappa)$ is positive and bounded away
from zero for $g \in \supp(\phi)$.
\end{lem}

\pf
There exists $r_0 \ge 1$ such that, for all $k=1,\dots,|W|$,
whenever $g \in \supp(\phi)$, $\zeta \in T^*X$,
$\|z(g,\zeta)\|_k \le 3D$ we have:
$$
\|\zeta\| \le 1  \text{ implies } \|\xi(g,\zeta)\|_k \le r_0
\quad \text{and} \quad
\|\xi(g,\zeta)\|_k \le 1 \text{ implies } \|\zeta\| \le r_0.
$$
Hence, under the same conditions, we also have:
$$
\|\zeta\| \le s \text{ implies }  \|\xi(g,\zeta)\|_k \le s r_0
\quad \text{and} \quad
\|\xi(g,\zeta)\|_k \le s  \text{ implies }  \|\zeta\| \le s r_0,
$$
for all $s>0$.

It follows by induction on $k$ that, whenever
$g \in \supp(\phi)$, $\zeta \in T^*X$, $\|\zeta\| \le s$,
we have
$\|(\Theta_t^k \circ \dots \circ \Theta_t^1)(g,\zeta)\|
\le s(r_0)^{2k}$.
Indeed, suppose that the statement holds for $k=l$.
If
$\|z((\Theta_t^l \circ \dots \circ \Theta_t^1)(g,\zeta))\|_{l+1} \ge 2D$,
then
$$
(\Theta_t^{l+1} \circ \Theta_t^l \circ \dots \circ \Theta_t^1)(g,\zeta)
=(\Theta_t^l \circ \dots \circ \Theta_t^1)(g,\zeta),
$$
and so 
\begin{multline*}
\|(\Theta_t^{l+1} \circ \Theta_t^l \circ \dots \circ \Theta_t^1)(g,\zeta)\|=
\|(\Theta_t^l \circ \dots \circ \Theta_t^1)(g,\zeta)\|   \\
\le s (r_0)^{2l} \le s (r_0)^{2(l+1)}.
\end{multline*}
Or else 
$\|z((\Theta_t^l \circ \dots \circ \Theta_t^1)(g,\zeta))\|_{l+1} \le 2D$,
in which case
$$
\|z((\Theta_t^{l+1} \circ \Theta_t^l \circ \dots \circ \Theta_t^1)
(g,\zeta))\|_{l+1} \le 3D,
$$
and so
\begin{multline*}
\|(\Theta_t^{l+1} \circ \Theta_t^l \circ \dots \circ \Theta_t^1)(g,\zeta)\| \\
\le r_0
\|\xi((\Theta_t^{l+1} \circ \Theta_t^l \circ \dots \circ \Theta_t^1)
(g,\zeta))\|_{l+1} \\
= r_0 \|\xi((\Theta_t^l \circ \dots \circ \Theta_t^1)(g,\zeta))\|_{l+1}  \\
\le r_0^2
\|(\Theta_t^l \circ \dots \circ \Theta_t^1)(g,\zeta)\|
\le s (r_0)^{2(l+1)},
\end{multline*}
which proves the induction step.

\separate

Let $\tilde r_0 = (r_0)^{1-2|W|}$.
Then, for all $k=1,\dots,|W|-1$,
whenever $g \in \supp(\phi)$, $\zeta \in T^*X$,
$\|\zeta\| \le \tilde r_0/t$, we have:
\begin{multline*}
\delta
\bigl( \|z((\Theta_t^{k-1} \circ \dots \circ \Theta_t^1)(g,\zeta))\|_k \bigr)
\gamma \bigl(
t\|\xi((\Theta_t^{k-1} \circ\dots\circ \Theta_t^1)(g,\zeta))\|_k \bigr)  \\
=\delta
\bigl( \|z((\Theta_t^{k-1} \circ \dots \circ \Theta_t^1)(g,\zeta))\|_k \bigr).
\end{multline*}
Write
\begin{multline*}
\langle g, \mu(\zeta) \rangle -
\Theta_t^* \langle g, \mu(\zeta) \rangle
= \bigl( \langle g, \mu(\zeta) \rangle -
\langle g, \mu(\Theta_t^1(g,\zeta)) \rangle \bigr)   \\
+ \dots + \bigl(
\langle g, \mu((\Theta_t^{k-1} \circ\dots\circ \Theta_t^1)(g,\zeta)) \rangle -
\langle g, \mu((\Theta_t^k \circ \dots \circ \Theta_t^1)(g,\zeta)) \rangle
\bigr)  \\
+ \dots + \bigl( \langle g,
\mu((\Theta_t^{|W|-1} \circ \dots \circ \Theta_t^1)(g,\zeta)) \rangle -
\langle g,\mu((\Theta_t^{|W|} \circ \dots \circ \Theta_t^1)(g,\zeta)) \rangle
\bigr).
\end{multline*}
Then let $(g,\zeta_k)= (\Theta_t^{k-1} \circ \dots \circ \Theta_t^1)(g,\zeta)$
and notice that, in the coordinate system $\tilde \psi_k$,
\begin{multline*}
\langle g, \mu(\zeta_k) \rangle -
\langle g, \mu(\Theta_t^k(g,\zeta_k)) \rangle
= t \epsilon \delta(\|z\|_k(g,\zeta_k))\gamma(t\|\xi(g,\zeta_k)\|_k) \\
\cdot \bigl( |\alpha_{x_k,1}(g)| \xi_1(g,\zeta_k)\bar \xi_1(g,\zeta_k) + \dots
+ |\alpha_{x_k,n}(g)| \xi_n(g,\zeta_k) \bar \xi_n(g,\zeta_k) \bigr)    \\
= t \epsilon \delta(\|z\|_k(g,\zeta_k))
\cdot \bigl( |\alpha_{x_k,1}(g)| |\xi_1(g,\zeta_k)|^2 + \dots
+ |\alpha_{x_k,n}(g)| |\xi_n(g,\zeta_k)|^2 \bigr)   \\
= t \|\zeta\|^2 \epsilon \delta(\|z\|_k(g,\zeta_k))
\cdot \frac{|\alpha_{x_k,1(g)}| |\xi_1(g,\zeta_k)|^2 + \dots
+ |\alpha_{x_k,n}(g)| |\xi_n(g,\zeta_k)|^2}{\|\zeta\|^2}.
\end{multline*}

Hence, when $t'\|\zeta\| \le \tilde r_0$, we get
$$
\kappa(g,\zeta, t') =
\langle g, \mu(\zeta) \rangle - \Theta_{t'}^* \langle g, \mu(\zeta) \rangle
= t'\|\zeta\|^2 \cdot
\tilde \kappa \bigl( g, \frac{\zeta}{\|\zeta\|}, t' \bigr)
$$
for some smooth bounded function
$\tilde \kappa (g, v, t')$
defined on $\Omega \times \{\zeta \in T^*X;\, \|\zeta\|=1 \} \times [0,1]$.
It is clear that $Re(\tilde \kappa)$ is positive.
Because the set
$\supp(\phi) \times \{\zeta \in T^*X;\, \|\zeta\|=1 \} \times [0,1]$
is compact, $Re(\tilde \kappa)$ is bounded away from zero on this set.
\qed

Thus in the region defined by $t'\nu \le \tilde r_0$ the function in
question becomes
$$
e^{-\kappa(g,\zeta, t')} =
e^{- t' \|\zeta\|^2 \cdot
\tilde \kappa \bigl( g, \frac{\zeta}{\|\zeta\|}, t' \bigr)}
$$
for some smooth bounded function $\tilde \kappa$
whose real part $Re(\tilde \kappa)$ is positive and
bounded away from zero.
Because the set
$\supp(\phi) \times \{\zeta \in T^*X;\, \|\zeta\|=1 \} \times [0,1]$
is compact, all the partial derivatives of $\tilde \kappa$ with respect
to the $g$ variable can be bounded on this set.
Hence it follows that all the partial derivatives of
$e^{-t'\nu ^2 \cdot \tilde \kappa}$
with respect to $g$ variable are bounded on the region
$g \in \supp(\phi)$, $t'\nu \le \tilde r_0$.

\separate

Finally, suppose that $t'\nu \ge \tilde r_0$.
We already know that the function $\frac 1{\|\zeta\|} \kappa$
can be written independently of the variable $\nu$.
Lemma \ref{scaling} implies that the function
$\frac 1{\|\zeta\|} \kappa$
becomes constant with respect to $t'\nu$ when $t'\nu$ is
larger than some constant $\tilde R_0$.
Thus all the partial derivatives of $\frac 1{\|\zeta\|} \kappa$
are bounded.

On the other hand, Lemma \ref{realpartlemma} implies
that $\frac 1{\|\zeta\|} Re(\kappa)$ is positive and bounded away
from zero for $g \in \supp(\phi)$, $t'\nu \ge \tilde r_0$.
Hence boundedness of the partial derivatives of
$\frac 1{\|\zeta\|} \kappa$ implies that the partial derivatives of
$e^{-\|\zeta\| \cdot \frac 1{\|\zeta\|} \kappa}$
with respect to $g$ variable are bounded on the region
$g \in \supp(\phi)$, $t'\nu \ge \tilde r_0$.

This shows that all the partial
derivatives of all orders of $e^{-\kappa(g, \zeta, t')}$
with respect to the $g$ variable can be bounded independently
of $\zeta$ and $t'$ on the set
$\supp(\phi) \times \{\zeta \in T^*X;\, \|\zeta\|>0\} \times [0,t]$.
Hence we proved Lemma \ref{slanting}.
\qed

\separate

\appendix

\begin{section}
{Equality of Coefficients}  \label{appendix}
\end{section}

In this appendix we will show that the coefficients $m_k$'s
coincide with coefficients $d_{g_0,x_k}$'s in \cite{SchV2}.
Recall that
\begin{multline*}
m_k = \chi \bigl( (Rj_{B^0_k \hookrightarrow X})_* \circ
(j_{B^0_k \hookrightarrow B_k})^! ({\cal F}|_{B_k})_{x_k} \bigr)  \\
= \chi \bigl( (j_{\{x_k\} \hookrightarrow O_k})^* \circ
(Rj_{B^0_k \hookrightarrow O_k})_* \circ
(j_{B^0_k \hookrightarrow O_k})^! ({\cal F}|_{O_k}) \bigr)  \\
= \chi \bigl( R\Gamma_{B^0_k} ({\cal F}|_{O_k})_{x_k} \bigr)
= \chi \bigl( R\Gamma_{\{x_k\}} ({\cal F}|_{O_k})_{x_k} \bigr)
+ \chi \bigl( R\Gamma_{B^0_k \setminus \{x_k\}} ({\cal F}|_{O_k})_{x_k} \bigr),
\end{multline*}
and formula (5.25b) from \cite{SchV2} says that
$$
d_{g_0,x_k} = \chi \bigl( {\cal H}^*_{O_k} (\BB D {\cal F})_{x_k} \bigr)
= \chi \bigl( (j_{x_k \hookrightarrow O_k})^* \circ 
(j_{O_k \hookrightarrow X})^! (\BB D {\cal F}) \bigr),
$$
where $\BB D {\cal F}$ is the Verdier dual of ${\cal F}$.
It follows from Proposition 3.1.13 and Example 3.4.5 of \cite{KaScha}
that
$$
(j_{O_k \hookrightarrow X})^! (\BB D {\cal F}) \simeq
\BB D_{O_k} ({\cal F}|_{O_k}).
$$
Hence by Proposition 3.4.3 of \cite{KaScha}
\begin{multline*}
d_{g_0,x_k} = \chi \bigl( \BB D_{O_k} ({\cal F}|_{O_k})_{x_k} \bigr) =
\chi \bigl( R\operatorname{Hom}
(R\Gamma_{\{x_k\}} (O_k, {\cal F}|_{O_k}), \BB C) \bigr)  \\
= \chi \bigl( R\Gamma_{\{x_k\}} (O_k, {\cal F}|_{O_k}) \bigr)
= \chi \bigl( R\Gamma_{\{x_k\}} ({\cal F}|_{O_k})_{x_k} \bigr).
\end{multline*}
Thus to show $m_k = d_{g_0,x_k}$ it is enough to show that
\begin{multline}  \label{eulerchar}
\chi \bigl( R\Gamma_{B^0_k \setminus \{x_k\}} ({\cal F}|_{O_k})_{x_k} \bigr) \\
= \chi \bigl( (Rj_{ B^0_k \setminus \{x_k\} \hookrightarrow O_k})_* \circ
(j_{ B^0_k \setminus \{x_k\} \hookrightarrow O_k})^!
({\cal F}|_{O_k})_{x_k} \bigr) = 0.
\end{multline}
Let
$$
B_{\epsilon} =
\psi_{g_0,k} \bigl(\{ n^0 \in \g n_k^0;\, \|n^0\|_k < \epsilon \}\bigr).
$$
Then
\begin{multline*}
(Rj_{ B^0_k \setminus \{x_k\} \hookrightarrow O_k})_* \circ
(j_{ B^0_k \setminus \{x_k\} \hookrightarrow O_k})^!
({\cal F}|_{O_k})_{x_k}   \\
= \varinjlim_{B_{\epsilon}} R\Gamma_{B_{\epsilon} \setminus \{x_k\}}
\bigl( B_{\epsilon} \setminus \{x_k\},({\cal F}|_{O_k}) \bigr).
\end{multline*}
The sheaf ${\cal F}$ being $G_{\BB R}$-equivariant and the inclusion (\ref{r})
together imply that, whenever $\epsilon <r$, each cohomology sheaf
${\cal H}^l \bigl( (j_{ B_{\epsilon} \setminus \{x_k\} \hookrightarrow O_k})^!
({\cal F}|_{O_k}) \bigr)$
is a multiple of the constant sheaf on $B_{\epsilon} \setminus \{x_k\}$.
Thus, because the Euler characteristic of $B_{\epsilon} \setminus \{x_k\}$
is zero,
$$
\chi R\Gamma_{B_{\epsilon} \setminus \{x_k\}}
\bigl( B_{\epsilon} \setminus \{x_k\},({\cal F}|_{O_k}) \bigr) =0,
$$
which proves (\ref{eulerchar}).
This finishes our proof that the coefficients $m_k$'s and
$d_{g_0,x_k}$'s are equal.

\separate

\end{document}